\documentclass[12pt]{article}
\usepackage{tikz}
\usepackage{graphicx}
\usepackage{multirow}
\usepackage{amsmath}
\usepackage{cases}
\usepackage{indentfirst}
\usepackage{mathrsfs}
\usepackage{amsfonts,amssymb}
\usepackage{amsthm} 
\usepackage[all]{xy}

\newtheorem{lem}{Lemma}[section]%
\newtheorem{thm}[lem]{Theorem}%
\newtheorem*{thm*}{Theorem}
\newtheorem{defi}[lem]{Definition}%
\newtheorem{cor}[lem]{Corollary}%
\newtheorem{exam}[lem]{Example}%
\newtheorem{rem}[lem]{Remark}%
\newenvironment{pf}{\medskip\noindent{Proof.\hspace{0.2cm}}}{\hfill \qed \newline \medskip}

\def\a{\alpha}   \def\d{\delta} 
    \def\k{\kappa}

\def\nd{\mathrel{\bigm|\kern-.7em/}}

\begin{document}
	
	\title	{\Large  Some Results on Cohomological Properties of $p$-groups and Non-inner Automorphisms with order $p$ for Non-abelian Finite $p$-groups }	
	
	\author{Wei Xu}
	
	\maketitle
	\begin{abstract}
		In this paper, we prove that any finite $p$-group of order at least $p^{2}$ admitting non-inner automorphism of order $p$. 
		
		We had originally planned to revise our arXiv preprint after receiving the journal’s referee reports. On November 22, an AI-assisted paper that  claim to have solved the conjecture appeared online, so we now have to update the arXiv version earlier than planned. We did not use AI in the entire paper.
	\end{abstract}
	\section{Preliminaries}
	In 1966, Wolfgang Gaschütz had proved that any finite $p$-group of order at least $p^{2}$ admitting outer automorphism of order $p$ in \cite{gas}. A natural problem was proposed by Ya.G.Berkovich: 	
	\textbf{	whether any finite $p$-group with order at least $p^{2}$ has non-inner automorphism of order $p$.} The above statement holds true for abelian $p$-groups with order at least $p^{2}$. Naturally, it suffices to solve the above problem for the non-interactive case, namely
	\textbf{	whether any non-abelian finite $p$-group has a non-inner automorphism of order $p$.} It was writed in problem 4.13 of \cite{kou}.\\
	
	Let $G$ be a non-abelian finite $p$-group. Marian Deaconescu and Gheorghe Silberberg had proved that any finite $p$-group $G$ with $C_{G}(\Phi(G))\nleqslant\Phi(G)$ admitting a non-inner automorphism of order $p$ in \cite{marian}. Alireza Abdollahi had proved that powerful $p$-group has a non-inner automorphism of order $p$ in \cite{alireza1}. 
	
	In this paper, we show that the answer to the above problem is positive. In section 2, Theorem \ref{5.5} is proved. From this, only discuss $\mathrm{H}^{1}(G/N,\Omega_{1}(Z(N)))\leq\mathbb{F}_{p}^{n}$ where $n=d(Z(G))$. In section 3, some definitions  are given and some key lemmas are presented. In section 4, some cohomological properties of $n-G$-module are verified. In section 5, by using cohomological properties proved in previous section, we get the following results:
	\begin{thm*}
		Suppose $N$ is a special subgroup of $G$. Then one of following holds:
		
		{\rm(1)}\; $G$ has a non-inner automorphism of order $p$;
		
		{\rm(2)}\; Exists a special subgroup $M$ of $G$ such that $|M|< |N|$;
		
		{\rm(3)}\; Exists a special subgroup $M$ of $G$ such that $|M|=|N|,|\mathrm{I}(C_{G}(N))|<|\mathrm{I}(C_{G}(M))|$.
	\end{thm*} 
	From the above result, the following Theorem can be proved.
	
	\begin{thm*}
		Any finite $p$-group of order at least $p^{2}$ has a non-inner automorphism of order $p$.
	\end{thm*}
	
	In this paper, all groups are finite $p$-groups unless otherwise stated. For a group $H$, $H\leq \mathrm{C}_{p}^{m}$ means $H\cong\mathrm{C}_{p}^{t}$ where $t\leq n$. The mean of $H\geq \mathrm{C}_{p}^{m}$ is similar. The following notations are standard in group theory.	$ \Phi(G)$: The intersection of all maximal subgroups of $G$;	$d(G)$: $|G/\Phi(G)|=p^{d(G)}$;	$\Omega_{1}(G)=\langle x\in G\mid x^{p}=1\rangle$;	$Z(G)=\{g\in G\mid h^{g}=h\;for\;all\;h\in G\}$.
	\section{Result on Non-inner Automorphism of Order $p$}
	
	\begin{defi}
		Let $G$ be a non-abelian finite $p$-group and $A$ be a subgroup of $G$. Denote
		$$\mathrm{I}(A)=\{g\in A\mid g^{p}\in Z(G)\}.$$ 
	\end{defi}
	Obviously, $ \mathrm{I}(A)$ is a group when $A$ is abelian.	
	In this paper, use the description of $ Z^{n}(G,M)$ and $ B^{n}(G,M)$ in Definition 16.12 in \cite{hu}.
	Let $N$ be a normal subgroup of $G$ with $N\leq\Phi(G)$. For $\tau\in Z^{1}(G/N,\Omega_{1}(Z(N)))$, set $\psi(g)=g\tau(gN)$ for all $g\in N$. From Lemma 2.8 in \cite{rus}, $\psi$ is an automorphism of order $p$. Say $\psi$ is automorphism induced by $\tau$. In fact, there is a simple Lemma about order of automorphism.
	
	\begin{lem}\label{lp}
		Let $N$ be a normal subgroup of $G$ with $N\leq\Phi(G)$. And $N_{1}$ is a normal subgroup of $G$ that $\Omega_{1}(Z(N))\leq N_{1}<N$. Set $\psi(g)=g\tau(gN_{1})$ for all $g\in G$ where $\tau\in Z^{1}(G/N_{1},\Omega_{1}(Z(N)))$. Then $\psi$ is order $p$.
	\end{lem}
	\begin{pf}
		Claim $\psi^{n}(g)=g(\tau(gN_{1}))^{n}$. By induction on $n$, suppose  $\psi^{t}(g)=g(\tau(gN_{1}))^{t}$. Then
		$$\psi^{t+1}(g)=\psi(\psi^{t}(g))=\psi(g(\tau(gN_{1}))^{t})=g(\tau(gN_{1}))^{t}(\tau(g(\tau(gN_{1}))^{t}N_{1})).$$ Since $\tau(gN_{1})\in N_{1}$ for all $g\in G$, thus $\psi^{t+1}(g)=g(\tau(g))^{t+1}$. The assertion is proved.
		
		Since $\tau\in Z^{1}(G/N_{1},\Omega_{1}(Z(N)))$, $\psi$ is order $p$.
		
	\end{pf}
	\begin{lem}\label{ij}
		Let $G$ be a non-abelian $p$-group, and let $ N$ be a  normal subgroup of $G$; $n=d(Z(G))$. If  $[N,N]\leq Z(G)\leq N$ and $\Omega_{1}(N)$ is abelian, then $$\mathrm{I}(N)/Z(G)\Omega_{1}(N)\leq \mathrm{C}_{p}^{n}.$$
		
	\end{lem}
	\begin{pf}
		Set 
		\begin{align*}
			\Delta: \mathrm{I}(N)\rightarrow Z(G)\\
			x\mapsto x^{p}.
		\end{align*}
		Denote $H=\{x\in \mathrm{I}(N)|x^{p}\in\mho_{1}(Z(G))\}$. Then $\Delta$ induce a new homomorphism 
		\begin{align*}
			\Delta_{1}: \mathrm{I}(N)/H&\rightarrow Z(G)/\mho_{1}(Z(G))\\
			xH&\mapsto x^{p}\mho_{1}(Z(G)).
		\end{align*}
		Notice $\Delta_{1}$ is injective.  Since $Z(G)/\mho_{1}(Z(G))\cong\mathrm{C}_{p}^{n} $ and $ H=\Omega_{1}(N)Z(G)$, $$\mathrm{I}(N)/\Omega_{1}(N)Z(G)\leq \mathrm{C}_{p}^{n}.$$ 
		
	\end{pf}
	
	\begin{lem}\label{ddd}
		Let $G$ be a non-abelian $p$-group and $n=d(Z(G))$. And let $N$ be a normal subgroup of $G$ such that 
		\begin{align*}
			&C_{G}(N)/Z(N)\;is\;cyclic,\; and\\
			&\mathrm{I}(C_{G}(N))\leq N\leq NC_{G}(N)\leq \Phi(G).
		\end{align*}
		Then $\Omega_{1}(Z(NC_{G}(N)))\leq N$. 
		
		Denote $W=\Omega_{1}(Z(NC_{G}(N)))$. 
		If all automorphisms induced by derivations in $Z^{1}(G/NC_{G}(N),W) $ are inner, then $\mathrm{H}^{1}(G/NC_{G}(N),W)\cong\mathrm{I}(N)/Z(G)\Omega_{1}(N) $.
	\end{lem}
	\begin{pf} Set $$ \mathrm{Z}=\{x\in G\mid g^{x}=g\tau(gNC_{G}(N))\;for\;all\;g\in G \;and\;some\;\tau\in Z^{1}(G/NC_{G}(N),W)\}$$
		and $$ \mathrm{B}=\{x\in G\mid g^{x}=g\tau(gNC_{G}(N))\;for\;all\;g\in G \;and\;some\;\tau\in B^{1}(G/NC_{G}(N),W)\}.$$
		Since all automorphisms induced by derivations in $Z^{1}(G/NC_{G}(N),W) $ fix $NC_{G}(N)$, $\mathrm{Z}\leq C_{G}(NC_{G}(N))\leq C_{G}(N)$. From all automorphisms induced by derivations in $Z^{1}(G/NC_{G}(N),W) $ are order $p$, $\mathrm{Z}\leq \mathrm{I}(C_{G}(N))$. Set $\d_{x}(gNC_{G}(N))=g^{-1}g^{x}$ for all $g\in G$ where $x\in \mathrm{I}(C_{G}(N))$. Then $\d_{x}\in Z^{1}(G/NC_{G}(N),W)$. Thus $\mathrm{Z}= \mathrm{I}(C_{G}(N))$. It's easy to prove $\mathrm{B}=Z(G)\Omega_{1}(N)$. Set
		\begin{align*}
			\rho:\mathrm{Z}/\mathrm{B}&\rightarrow Z^{1}(G/NC_{G}(N),W)/B^{1}(G/NC_{G}(N),W)\\
			x\mathrm{B}&\mapsto \d_{x}+B^{1}(G/NC_{G}(N),W)
		\end{align*}
		where $\d_{x}(gNC_{G}(N))=g^{-1}g^{x}$ for all $g\in G$.
		From the definition of $\mathrm{Z} $ and $ \mathrm{B}$, $\rho$ is bijective. Thus $\mathrm{Z}/\mathrm{B}\cong\mathrm{H}^{1}(G/NC_{G}(N),W).$
		
	\end{pf}
	
	\begin{thm}\label{5.5}
		Let $G$ be a non-abelian $p$-group with $n=d(Z(G))$. And $N$ is a normal subgroup of $G$ that 
		\begin{align*}
			&C_{G}(N)/Z(N)\;is\;cyclic,\; and\\
			&\mathrm{I}(C_{G}(N))\leq N\leq NC_{G}(N)\leq \Phi(G).
		\end{align*}
		Then $\Omega_{1}(Z(NC_{G}(N)))\leq N$.  If\; $\mathrm{H}^{1}(G/NC_{G}(N),\Omega_{1}(Z(N)))\geq\mathbb{F}_{p}^{n+1}$, then $G$ has a non-inner automorphism of order $p$.

	\end{thm}
	\begin{pf}
		Denote $W=\Omega_{1}(Z(NC_{G}(N)))$ and $A=NC_{G}(N)$. Since $C_{G}(N)/Z(N) $ is cyclic, $C_{G}(N)$ is abelian. Notice $ C_{G}(NC_{G}(N))\leq C_{G}(N)$. Then $C_{G}(N)=Z(NC_{G}(N))=C_{G}(NC_{G}(N))$.
		Suppose all automorphisms induced by derivations in $Z^{1}(G/A,W)$ are inner.
		
		If $\mathrm{H}^{1}(G/A,W)\geq\mathbb{F}_{p}^{n+1}$, from Lemma \ref{ddd},  $\mathrm{I}(C_{G}(A))/WZ(G)\geq\mathrm{C}_{p}^{n+1}$.  It's contradictory to Lemma \ref{ij}.
		
		The proof is completed.
	\end{pf}
	
	\section{Some Basic Lemmas and Labels}

	Set $\tau\in Z^{2}(G,\mathrm{C}_{p})$. Denote $_{\tau}G.\mathrm{C}_{p}$ is the extension of  $\mathrm{C}_{p}$ by $G$ constructed from $\tau$ in chapter 11 of \cite{rb}. And $M$ is a $G$-module. For a subgroup $N$ of $G$, set $M^{N}=\{a\in M|a^{g}=a\;for\;all\;g\in N\}$. And $\mathbb{F}_{p}(G)$ is group algebra of $G$ over $\mathbb{F}_{p}$.
	
	From Theorem \ref{5.5}, only discuss $\mathrm{H}^{1}(G/N,\Omega_{1}(Z(N)))\leq\mathbb{F}_{p}^{n}$ where $n=d(Z(G))$. Then give the following definition. 
	\begin{defi}
		Let $G$ be a finite $p$-group. Say a $\mathbb{F}_{p}(G)$-module $ A$ is $n-G$-module if $A^{G}\cong \mathbb{F}_{p}^{n}$ and $\mathrm{H}^{1}(G,A)\leq \mathbb{F}_{p}^{n}$. Specially, when $\mathrm{H}^{1}(G,A)\cong\mathbb{F}_{p}^{n}$, say $A$ is exactly $n-G$-module.
	\end{defi}
	
	From Theorem \ref{5.5},	a key problem on dealing with non-inner automorphism of order $p$ on non-abelian finite $p$-group is that 
	\begin{center}
		\textbf{	whether $\mathrm{H}^{1}(_{\tau}G.\mathrm{C}_{p},A)\ncong \mathrm{H}^{1}(G,A)$ for $\tau\in Z^{2}(G,\mathrm{C}_{p})$}
	\end{center} when $A$ is a $n-G$-module? In this section 4, we will prove some results about the problem. Specificly, we will embed $ A$ in direct sum of $n$ copies of $\mathbb{F}_{p}(G)$. And disscuss some properties on direct sum of $n$ copies of $\mathbb{F}_{p}(G)$. 
	
	In the following section, preparing some basic Lemmas for discussion in section 4.

	\begin{defi}
		Assume $A$ be a $\mathbb{F}_{p}(G)$-module. Let 
		$$\mathrm{J}_{G}(A)=\{\sum\limits_{i=1}\limits^{m}x_{i}^{g_{i}}|x_{i}\in A, g_{i}\in\mathrm{J}(\mathbb{F}_{p}(G)) \;for\;all\;1\leq i\leq n\}$$ where $\mathrm{J}(\mathbb{F}_{p}(G)) $ is Jacobson radical of $\mathbb{F}_{p}(G) $.
	\end{defi}
	Sometimes, $G$-module $A$ is both right $G$ module and left $G$-module. And $\mathrm{J}_{G}(A)$ can be different in two cases. So before calculating $\mathrm{J}_{G}(A)$, it's necessary to determine the state of $A$.

	\begin{defi}
		Assume $A$ is a $\mathbb{F}_{p}(G)$ module. And $S=\{x_{1},\cdots,x_{n}\}$ is a subset of $A$. Say $S$ is a $G$ minimal generator set of $A$ if $A$ is a $\mathbb{F}_{p}(G)$ module generated by $S$ and $\mathbb{F}_{p}(G)$-module generated by aribitrary proper subset of $S$ is a proper $\mathbb{F}_{p}(G)$-submodule of $A$.
	\end{defi}
	
	\begin{lem}\label{1.2}
		Let $A$ be a $\mathbb{F}_{p}(G)$-module. And $A/\mathrm{J}_{G}(A)\cong\mathbb{F}_{p}^{n}$. Assume $A$ has a minimal generator set $ S$. Then $|S|=n$.
	\end{lem}
	\begin{pf}
		Suppose $S=\{x_{1},\cdots,x_{m}\}$. Set $B_{i}=\{\sum\limits_{i=1}\limits^{m}x_{i}^{g_{i}}|g_{i}\in \mathbb{F}_{p}(G)\}$ and $A_{i}=\{\sum\limits_{i=1}\limits^{m}x_{i}^{g_{i}}|g_{i}\in \mathrm{J}(\mathbb{F}_{p}(G))\}$. Notice $\mathbb{F}_{p}(G)$ is a local ring. Then $A_{i}$ is unique maximal $\mathbb{F}_{p}(G)$ submodule of $B_{i}$. Since $S$ is a minimal generator set of $A$, $ A/(A_{1}+\cdots+A_{m})$ is a direct sum of $m$ simple $ \mathbb{F}_{p}(G)$-modules. And any simple $ \mathbb{F}_{p}(G)$-module is isomorphic to $\mathbb{F}_{p}$ as trivial  $ \mathbb{F}_{p}(G)$- module. Then 
		$$ A/(A_{1}+\cdots+A_{m})\cong\mathbb{F}_{p}^{m}$$ as $ \mathbb{F}_{p}(G)$-module. Notice $G$ acts on $A/(A_{1}+\cdots+A_{m})$ trivially. Then $\mathrm{J}_{G}(A)\leq A_{1}+\cdots+A_{n}$. And $m\leq n$.
		
		Notice $G$ acts $ A/\mathrm{J}_{G}(A)$ trivial. Then $m\geq n$. Now $m=n$.
	\end{pf}
	
	Now, give a new definition.
	
	\begin{defi}\label{up}
		Assume $A$ is a $\mathbb{F}_{p}(G)$-module. Let 
		$$ d_{G}(A)$$
		is positive integer that $A/\mathrm{J}_{G}(A)\cong\mathbb{F}_{p}^{d_{G}(A)}$.
	\end{defi}
	
	The following lemmas will be used in Section 2.
	
	\begin{lem}\label{cc}
		Let $A$ be a $\mathbb{F}_{p}(G)$ module. Set $A_{1}$ is a  $\mathbb{F}_{p}(G)$ proper submodule of $A$ that $d_{G}(A)=m$. Assume $d_{G}(A/A_{1})=s$. Then $d_{G}(A_{1})\leq s+m$.
	\end{lem}
	\begin{pf}
		Denote $\{x_{1},\cdots,x_{m}|x_{i}\in A_{1}\}$ as a minimal generator set of $ A_{1}$. And set $y_{1},\cdots,y_{s}\in A$ that $ A/A_{1}$ is generated by $y_{1}+A_{1},\cdots,y_{s}+A_{1}$ as  $\mathbb{F}_{p}(G)$-module. Then $A$ is generated by $x_{1},\cdots,x_{m},y_{1},\cdots,y_{s}$ as $\mathbb{F}_{p}(G)$-module. Thus $d_{G}(A_{1})\leq s+m$.
	\end{pf}
	
	When $\mathbb{F}_{p}$ is recorded as $G$-module, $G$ acts on $\mathbb{F}_{p}$ trivially by automorphism group of additive group $\mathbb{F}_{p}$ is isomorphic to $\mathrm{C}_{p-1}$. Next, $\bigoplus\limits^{n}\mathbb{F}_{p}(G)$ is direct sum of $n$ copies of $\mathbb{F}_{p}(G)$ as $\mathbb{F}_{p}(G)$-module.
	
	\begin{lem}\label{ut}
		Let $A$ be a $\mathbb{F}_{p}(G)$-module that $ A^{G}\cong\mathbb{F}_{p}^{n}$. Then $A$ can be embedded in the direct sum of $n$ copies of $\mathbb{F}_{p}(G)$ as $\mathbb{F}_{p}(G)$-module.
	\end{lem}
	\begin{pf}
		Notice $\mathbb{F}_{p}(G)$ is Frobenius algebra. From Proposition 1.6.2 in \cite{benson}, any $\mathbb{F}_{p}(G)$ free module generated by  finite elements is injective module. Since $A^{G}\cong (\bigoplus\limits^{n}\mathbb{F}_{p}(G))^{G}$, there exists $\mathbb{F}_{p}(G)$-module homomorphism  $\psi$ from $A$ to $\bigoplus\limits^{n}\mathbb{F}_{p}(G) $ such that $\psi(A^{G})=(\bigoplus\limits^{n}\mathbb{F}_{p}(G))^{G}$. From any simple $\mathbb{F}_{p}(G) $-module is isomorphic to $\mathbb{F}_{p}$, $|A_{1}\cap A^{G}|\geq p$ for any $\mathbb{F}_{p}$ submodule $A_{1}$ of $A$. Then $ker\psi=0$. Otherwise, $ker\psi\cap A^{G}\neq 0$. It's a contradiction. Thus $\psi$ is a embedded map.
		
		The proof is completed.
	\end{pf}

	\begin{cor}\label{1.00}
		Let $A$ be a $\mathbb{F}_{p}(G)$-module that $\mathrm{H}^{1}(G,A)=0$. Then $A\cong\bigoplus\limits^{n}\mathbb{F}_{p}(G)$ as $G$-module where $n$ is a positive integer.
	\end{cor}
	\begin{pf}
		Since all cohomological group of free module is trivial, only verify $A$ is free.
		Suppose $A^{G}\cong\mathbb{F}_{p}^{n}$. From Lemma \ref{ut}, $A$ can be embedded in $\bigoplus\limits^{n}\mathbb{F}_{p}(G)$ as as $G$-module where $n$ is a positive integer. Set $M$ is $G$-submodule of $\bigoplus\limits^{n}\mathbb{F}_{p}(G) $ that $M\cong A$ as $G$-module.
		
		If $A\cong\bigoplus\limits^{n}\mathbb{F}_{p}(G) $, from any simple $\mathbb{F}_{p}(G) $-module is isomorphic to $\mathbb{F}_{p}$, there exists $a\in \bigoplus\limits^{n}\mathbb{F}_{p}(G)\backslash M$ such that $a^{g-1}\in M$ for all $g\in G$. Set $\tau(g)=a^{g-1}$. Then $ \tau\in Z^{1}(G,M)\backslash B^{1}(G,M)$. It's contradictory to $\mathrm{H}^{1}(G,A)=0$.
		
		Now, the proof is completed.
	\end{pf}

	\begin{lem}\label{xxx}
		Let $A$ be a $\mathbb{F}_{p}(G)$-module that $d_{G}(A)=n$. Then $\mathrm{Hom}_{\mathbb{Z}}(A,\mathbb{F}_{p})^{G}\cong\mathbb{F}_{p}^{n}$
	\end{lem}
	\begin{pf}
		Set $\psi\in \mathrm{Hom}_{\mathbb{Z}}(A,\mathbb{F}_{p})^{G}$. Then $\psi(x^{g-1})=\,^{g-1}\psi(x)=0$ for all $x\in A,g\in G$. Thus $\phi(x)=0$ for all $\phi\in \mathrm{Hom}_{\mathbb{Z}}(A,\mathbb{F}_{p})^{G},x\in\mathrm{J}_{G}(A)$.
		Since $A/\mathrm{J}_{G}(A)\cong\mathbb{F}_{p}^{n}$, $\mathrm{Hom}_{\mathbb{Z}}(A,\mathbb{F}_{p})^{G}\cong\mathbb{F}_{p}^{n}$.
	\end{pf}

	\begin{lem}\label{1.22}
		Let $A$ be a $\mathbb{F}_{p}(G)$-module. Then $ \mathrm{Hom}_{\mathbb{Z}}(A,\mathbb{F}_{p})$ has $G$ minimal generator set of order $n$ if and only if $A^{G}\cong\mathbb{F}_{p}^{n}$ . 
	\end{lem}
	\begin{pf} Notice $A\cong\mathrm{Hom}_{\mathbb{Z}}(\mathrm{Hom}_{\mathbb{Z}}(A,\mathbb{F}_{p}),\mathbb{F}_{p})$ as $\mathbb{F}_{p}(G)$-module. When $ \mathrm{Hom}_{\mathbb{Z}}(A,\mathbb{F}_{p})$ has minimal generator set of order $n$, from Lemma \ref{xxx}, $ A^{G}\cong\mathbb{F}_{p}^{n}$.

		Suppose $A^{G}\cong\mathbb{F}_{p}^{n}$. From Lemma \ref{ut}, $A$ can be embeded in direct sum of $n$ copies of $\mathbb{F}_{p}(G)$ as $\mathbb{F}_{p}(G)$-module. 
		Then there is a epimorphism from $\mathrm{Hom}_{\mathbb{Z}}(\mathbb{F}_{p}(G)^{n},\mathbb{F}_{p}) $ to $\mathrm{Hom}_{\mathbb{Z}}(A,\mathbb{F}_{p}) $.
		Notice $\mathrm{Hom}_{\mathbb{Z}}(\mathbb{F}_{p}(G)^{n},\mathbb{F}_{p})\cong\mathbb{F}_{p}(G)^{n}$ as $ \mathbb{F}_{p}(G)$-module. Then $\mathrm{Hom}_{\mathbb{Z}}(A,\mathbb{F}_{p})$ has a minimal generator set of  order at least  $n$. Notice $^{\sum\limits_{g\in G}\lambda_{g}g}\psi(x)=\psi(x^{\sum\limits_{g\in G}\lambda_{g}g})=0$ for all $\sum\limits_{g\in G}\lambda_{g}g\in \mathrm{J}(\mathbb{F}_{p}(G)),\psi\in \mathrm{Hom}_{\mathbb{Z}}(A,\mathbb{F}_{p}),x\in A^{G}$. Then $$\mathrm{Hom}_{\mathbb{Z}}(A,\mathbb{F}_{p})/\mathrm{J}_{G}(\mathrm{Hom}_{\mathbb{Z}}(A,\mathbb{F}_{p}))\geq\mathbb{F}_{p}^{n}.$$ Thus $ \mathrm{Hom}_{\mathbb{Z}}(A,\mathbb{F}_{p})$ has $G$ minimal generator set of order $n$.
		
		The proof is completed.
	\end{pf}
	
	When $M$ is a $G$-module, denote $\pi$  is surjective homomorphism from $G.N$ to $G$. And set $m^{g}=m^{\pi(g)}$ for $g\in G.N$ and $m\in M$. Then $M$ is a $G.N$-module.

	\section{Some Results of Cohomological Group of $p$-groups}
	The newly added labels in this section will be only applied in this section.

	To get automorphism with order $p$, for $n-G$-module $A$, we hope get $\mathrm{H}^{1}(_{\tau}G.\mathrm{C}_{p},A)\ncong \mathrm{H}^{1}(G,A)$ for $\tau\in Z^{2}(G,\mathrm{C}_{p})\backslash B^{2}(G,\mathrm{C}_{p})$. And we will have the following Theorem.
	
	\begin{thm}\label{gg}
		Set $\tau\in Z^{2}(G,\mathrm{C}_{p})\backslash B^{2}(G,\mathrm{C}_{p})$. Let $Q$ be a $n-G$-module that $\mathrm{H}^{1}(G,Q)\cong\mathbb{F}_{p}^{m} $ where $m<n$.
		Then
		$\mathrm{H}^{1}(_{\tau}G.\mathrm{C}_{p},Q)\geq \mathrm{H}^{1}(G,Q)\oplus\mathbb{F}_{p}^{n-m}$.
	\end{thm}
	Since $Aut(\mathrm{C}_{p})\cong\mathrm{C}_{p-1}$, $G$ acts on $\mathrm{C}_{p}$ trivially.
	When $\mathrm{H}^{1}(G,A)\cong\mathbb{F}_{p}^{n}$, we have the following example.
	\begin{exam}\label{kk}
		Let $G\cong \mathrm{C}_{p}$. $M$ is Jacobson radical of $\mathbb{F}_{p}(G)$. Then $\mathrm{H}^{1}(G,M)\cong\mathbb{F}_{p} $. Set $ \tau\in Z^{2}(G,\mathrm{C}_{p})\backslash B^{2}(G,\mathrm{C}_{p})$. Then $_{\tau}G.\mathrm{C}_{p}$ is a cyclic $p$-group of $p^{2}$ order. And $ \mathrm{H}^{1}(_{\tau}G.\mathrm{C}_{p},M)\cong \mathrm{H}^{1}(G,M)$.
	\end{exam}
	
	Let $Q$ be a $n-G$-module. $P$ is a $G$ submodule of $Q$. Then $P$ can be naturally  recorded as a normal subgroup of $_{\tau}G.Q$.
	And we will prove the following Theorem.
	
	\begin{thm}\label{qq}
		Let $P$ be a $\mathbb{F}_{p}(G)$-module that $|P|=p^{2}$. And $P_{1}$ is a $G$-submodule of $P$ that $|P_{1}|=p$.	Set $\tau\in Z^{2}(G,P) $. 	Let $Q$ be a $n-G$-module. Then
		
		{\rm(1)} when $\mathrm{H}^{1}(G,Q)\cong\mathbb{F}_{p}^{m}$ that $ m<n$, $\mathrm{H}^{1}(_{\tau}G.P,Q)\geq \mathrm{H}^{1}(G,Q)\oplus \mathbb{F}_{p}^{2n-2m+1}$ for $m\geq 2$ and $\mathrm{H}^{1}(_{\tau}G.P,Q)\cong\mathbb{F}_{p}^{2n}$ for $m=0$.
		
		{\rm(2)} when $\mathrm{H}^{1}(G,Q)\cong\mathbb{F}_{p}^{n}$ and $\mathrm{H}^{1}(_{\tau}G.P/P_{1},Q)\cong\mathrm{H}^{1}(G,Q)$ for $ \tau\in Z^{2}(G,P)$, $ \mathrm{H}^{1}(_{\tau}G.P,Q)\cong \mathbb{F}_{p}^{2n}$.
	\end{thm}

	For a $n-G$-module $M$, $ M$ can be embedded in direct summand of $n$ copies of $\mathbb{F}_{p}(G)$. Then we will discuss the problem in submodule of direct summand of $n$ copies of $\mathbb{F}_{p}(G)$.
	\begin{defi}
		Let $ \prod\limits^{n}\mathbb{F}_{p}(G)$ be direct product of $n$ copies of $\mathbb{F}_{p}(G)$. 
		Multiplication and addition of $\prod\limits^{n}\mathbb{F}_{p}(G)$ are in each coordinate. Then $ \prod\limits^{n}\mathbb{F}_{p}(G)$ is an algebra over $\mathbb{F}_{p}$. 
		And \begin{align*}
			(x_{1},\cdots,x_{n})^{y}=(x_{1},\cdots,x_{n})(y,\cdots,y)=(x_{1}y,\cdots,x_{n}y)
		\end{align*}
		for all $x_{1},\cdots,x_{n},y\in\mathbb{F}_{p}(G)$. This give a right $\mathbb{F}_{p}(G)$-module structure on $\prod\limits^{n}\mathbb{F}_{p}(G)$. Left $\mathbb{F}_{p}(G)$-module structure on $\prod\limits^{n}\mathbb{F}_{p}(G)$ is similar. And $\prod\limits^{n}\mathbb{F}_{p}(G)$ is isomorphic to direct product of $n$ copies of $\mathbb{F}_{p}(G)$ as $\mathbb{F}_{p}(G)$-module.
		Set 
		\begin{align*}
			\Delta_{G}:\mathbb{F}_{p}(G)&\rightarrow \mathbb{F}_{p}\\
			\sum\limits_{g\in G}k_{g}g&\mapsto k_{1_{G}}
		\end{align*}
		and 
		\begin{align*}
			(x_{1},\cdots,x_{n})\cdot(y_{1},\cdots,y_{n})=\sum\limits_{i=1}\limits^{n}x_{i}y_{i}
		\end{align*}
		for all $(x_{1},\cdots,x_{n}),(y_{1},\cdots,y_{n})\in \prod\limits^{n}\mathbb{F}_{p}(G)$.
	\end{defi}
	
	\begin{lem}\label{00}
		Set $\mathscr{R}_{n,G}(\mathscr{L}_{n,G})$ is the set consist of all right(left) $G$-submodule of $\prod\limits^{n}\mathbb{F}_{p}(G)$.
		Set
		\begin{align*}
			\mathcal{L}_{G}:\mathscr{R}_{n,G}&\rightarrow \mathscr{L}_{n,G}\\
			Q&\mapsto \mathcal{L}_{G}(Q)
		\end{align*}
		where $\mathcal{L}_{G}(Q)=\{x\in \prod\limits^{n}\mathbb{F}_{p}(G)|x\cdot y=0\;for\;all\;y\in Q\}$. $\mathcal{L}_{G} $ is a bijective map. And $\mathcal{L}_{G}(Q)\cong \mathrm{Hom}_{\mathbb{Z}}(\prod\limits^{n}\mathbb{F}_{p}(G)/Q,\mathbb{F}_{p})$ as left $G$-module.
	\end{lem}
	\begin{pf}
		When $\prod\limits^{n}\mathbb{F}_{p}(G)$ is recorded as right $G$-submodule, $\mathrm{Hom}_{\mathbb{Z}}(\prod\limits^{n}\mathbb{F}_{p}(G),\mathbb{F}_{p})$ is a left $G$ module. For a right $G$-submodule $Q$ of $\prod\limits^{n}\mathbb{F}_{p}(G)$, $\mathrm{Hom}_{\mathbb{Z}}(\prod\limits^{n}\mathbb{F}_{p}(G)/Q,\mathbb{F}_{p})$ is a left $G$-submodule of $\mathrm{Hom}_{\mathbb{Z}}(\prod\limits^{n}\mathbb{F}_{p}(G),\mathbb{F}_{p})$.
		Set 	\begin{align*}
			\d:\prod\limits^{n}\mathbb{F}_{p}(G)&\rightarrow \mathrm{Hom}_{\mathbb{Z}}(\prod\limits^{n}\mathbb{F}_{p}(G),\mathbb{F}_{p})\\
			(x_{1},\cdots,x_{n})&\mapsto \d(x_{1},\cdots,x_{n})
		\end{align*}
		where $\d(x_{1},\cdots,x_{n})(y_{1},\cdots,y_{n})=\Delta(\sum\limits_{i=1}\limits^{n}x_{i}y_{i})$. It's easy to prove $\d$ is a left $G$- module isomorphism. Since $ Q$ is a right $G$-submodule of $\prod\limits^{n}\mathbb{F}_{p}(G)$, $\mathcal{L}_{G}(Q)=\{(x_{1},\cdots,x_{n})\in \prod\limits^{n}\mathbb{F}_{p}(G)|\Delta(\sum\limits_{i=1}\limits^{n}x_{i}y_{i})=0\;for\;all\;(y_{1},\cdots,y_{n})\in Q\}$. And $\mathcal{L}_{G}(Q)\cong \mathrm{Hom}_{\mathbb{Z}}(\prod\limits^{n}\mathbb{F}_{p}(G)/Q,\mathbb{F}_{p})$ as left $G$-module.
		Since $$ \mathrm{dim}_{\mathbb{F}_{p}}\mathrm{Hom}_{\mathbb{Z}}(\prod\limits^{n}\mathbb{F}_{p}(G)/Q,\mathbb{F}_{p})=\mathrm{dim}_{\mathbb{F}_{p}}\prod\limits^{n}\mathbb{F}_{p}(G)/Q,$$ $|\mathcal{L}_{G}(Q)||Q|=|\prod\limits^{n}\mathbb{F}_{p}(G)|$.
		
		Set
		\begin{align*}
			\mathcal{R}_{G}:\mathscr{L}_{n,G}&\rightarrow \mathscr{R}_{n,G}\\
			T&\mapsto \mathcal{R}_{G}(T)
		\end{align*}
		where $\mathcal{R}_{G}(T)=\{x\in \prod\limits^{n}\mathbb{F}_{p}(G)|y\cdot x=0\;for\;all\;y\in T\}$. Similar to the above disscussion, $|T|\cdot|\mathcal{R}_{G}(T)|=|\prod\limits^{n}\mathbb{F}_{p}(G)|$. For $Q\in\mathscr{R}_{n,G}$, $Q\leq\mathcal{R}_{G}\mathcal{L}_{G}(Q)$. Since $ |\prod\limits^{n}\mathbb{F}_{p}(G)|=|Q|\cdot|\mathcal{L}_{G}(Q)|=|\mathcal{L}_{G}(Q)|\cdot|\mathcal{R}_{G}\mathcal{L}_{G}(Q)|$, $Q=\mathcal{R}_{G}\mathcal{L}_{G}(Q) $. Then $\mathcal{R}_{G}\mathcal{L}_{G}$ is identity map on $ \mathscr{R}_{n,G}$. It's similar to prove that $\mathcal{L}_{G}\mathcal{R}_{G}$ is identity map on $ \mathscr{L}_{n,G}$. Thus $\mathcal{L}_{G}$ is a bijiective map.
	\end{pf}

	\begin{lem}\label{ww}
		For  a right(left) $G$-submodule $Q$ of $\prod\limits^{n}\mathbb{F}_{p}(G)$ that $Q^{G}\cong\mathbb{F}_{p}^{n}$, $\mathrm{H}^{1}(G,Q)=\mathbb{F}_{p}^{m}$ if and only if $\mathcal{L}_{G}(Q)(\mathcal{R}_{G}(Q))$ has $G$ minimal generator set of order $m$.
	\end{lem}
	\begin{pf}
		Since $\mathrm{H}^{1}(G,Q)=\mathbb{F}_{p}^{m}$, $ (\prod\limits^{n}\mathbb{F}_{p}(G)/Q)^{G}\cong\mathbb{F}_{p}^{m}$.
		From Lemma \ref{1.22}, $$\mathrm{Hom}_{\mathbb{Z}}(\prod\limits^{n}\mathbb{F}_{p}(G)/Q,\mathbb{F}_{p})$$ has a $G$ minimal generator set of order $m$. Since $\mathcal{L}_{G}(Q)\cong \mathrm{Hom}_{\mathbb{Z}}(\prod\limits^{n}\mathbb{F}_{p}(G)/Q,\mathbb{F}_{p})$, $\mathcal{L}_{G}(Q)$ has a $G$ minimal generator set of order $m$.  
		
		Since $\mathcal{L}_{G}(Q)$ has a minimal generator set of order $n$, $(\prod\limits^{n}\mathbb{F}_{p}(G)/Q)^{G}\cong\mathbb{F}_{p}^{m}$. Then $\mathrm{H}^{1}(G,Q)=\mathbb{F}_{p}^{m}$.
		
		Thep proof is completed.
	\end{pf}
	
	Let $P$ be a $\mathbb{F}_{p}(G)$-module that $|P|=p^{t}$.	Set $\tau\in Z^{2}(G,P)\backslash B^{2}(G,P)$.	Denote $\eta_{_{\tau},G}$ is surjective homomorphism from $_{\tau}G.P $  to $G$. Assume $ker\eta_{_{\tau},G}=\langle a_{1}\rangle\times\cdots\times\langle a_{t}\rangle$.

	Set
	\begin{align*}
		\Psi_{G}^{_{\tau}G.P}:\mathbb{F}_{p}(_{\tau}G.P)&\rightarrow \mathbb{F}_{p}(G)\\
		\sum\limits_{g\in _{\tau}G.P}k_{g}g&\mapsto\sum\limits_{g\in G}(\sum\limits_{\eta_{\tau,G}(h)=g}k_{h})g.
	\end{align*}
	$ \Psi_{G}^{_{\tau}G.P}$ is a algebra homomorphism.  Set
	\begin{align*}
		\Psi_{ _{\tau}G.P}^{G}:\mathbb{F}_{p}(G)&\rightarrow \mathbb{F}_{p}(_{\tau}G.P)\\
		\sum\limits_{g\in G}k_{g}g&\mapsto \sum\limits_{g\in G}k_{g}h_{g}(a_{1}-1)^{p-1}\cdots(a_{t}-1)^{p-1}
	\end{align*}
	where $h_{g}\in\mathbb{F}_{p}(_{\tau}G.P)$ that $\Psi_{G}^{_{\tau}G.P}(h_{g})=g$.
	For  $g\in ker\Psi_{G}^{_{\tau}G.P}$, $ g=g_{1}(a_{1}-1)+\cdots+g_{t}(a_{t}-1)$ where $g_{1},\cdots,g_{t}\in\mathbb{F}_{p}(_{\tau}G.P)$, then 
	$\Psi_{ _{\tau}G.P}^{G}$ is well defined.

	\begin{lem}\label{xo}
		Set 
		\begin{align*}
			_{n}\Psi_{G}^{_{\tau}G.P}:\prod\limits^{n}\mathbb{F}_{p}(_{\tau}G.P)&\rightarrow \prod\limits^{n}\mathbb{F}_{p}(G)\\
			(g_{1},\cdots,g_{n})&\mapsto(\Psi_{G}^{_{\tau}G.P}(g_{1}),\cdots,\Psi_{G}^{_{\tau}G.P}(g_{n})).
		\end{align*} and $ L=\{i\in\mathbb{N}| 1\leq p-1\}$ and $ U=\{i\in\mathbb{N}_{+}|i\leq n\}$. Let $L^{t}\times U$ be Cartesian product of $U$ and $t$ copies of $L$. Set  $e_{i_{1},\cdots,i_{t},l}$ is the element of $\prod\limits^{n}\mathbb{F}_{p}(_{\tau}G.P)$ that $l.th$ component is $ (a_{1}-1)^{i_{1}}\cdots(a_{t}-1)^{i_{t}}$ and others components are zero where $ 0\leq i_{1},\cdots,i_{t}\leq p-1,0\leq j\leq p-1,1\leq l\leq n$ and $(a_{i}-1)^{0}=1_{_{\tau}G.P}$ for all $1\leq i\leq n$. Let $ CR(G,\,_{\tau}G.P)$ be a set of coset representative of $_{\tau}G.P$ on $P$. Set $W=\{(x,\cdots,x)\in\prod\limits^{n}\mathbb{F}_{p}(_{\tau}G.P)|x\in CR(G,\,_{\tau}G.P)\}$. Set $W_{n}$ is $\mathbb{F}_{p}$ linear space that $W$ is a $\mathbb{F}_{p}$ basis of $W_{n}$.
		For $y\in ker _{n}\Psi_{G}^{_{\tau}G.P}$, there exists unique map $\lambda$ from $L^{n}\times U$ to $W_{n}$ such that 
		$$y=\sum\limits_{l=1}^{n}\sum\limits_{i_{1}+\cdots+i_{t}\geq1}\lambda(i_{1},\cdots,i_{t},l)e_{i_{1},\cdots,i_{t},l}.$$ Similarly, there exists unique map $\lambda_{1}$ from $L^{n}\times U$ to $W_{n}$ such that 
		$$y=\sum\limits_{l=1}^{n}\sum\limits_{i_{1}+\cdots+i_{t}\geq1}e_{i_{1},\cdots,i_{t},l}\lambda_{1}(i_{1},\cdots,i_{t},l).$$
	\end{lem}
	\begin{pf}
		Since mulitipulication of $\prod\limits^{n}\mathbb{F}_{p}(_{\tau}G.P)$ is mulitipulication on each component, only verify $n=1$. 
		
		Claim $y$ has a unique map $\k$ from $ CR(G,\, _{\tau}G.P) $ to $\mathrm{J}_{P}(\mathbb{F}_{p}(P))$ that $$ y=\sum\limits_{x\in CR(G,\, _{\tau}G.P)} x\k(x).$$ Since $ P$ is a normal subgroup of $ _{\tau}G.P$, there exists unique map $\tau$ from $CR(G,\, _{\tau}G.P)$ to $\mathbb{F}_{p}(P) $ such that $y=\sum\limits_{x\in CR(G,\, _{\tau}G.P)}x\tau(x)$. Since $y\in ker\Psi_{G}^{_{\tau}G.P}$, $\tau(x)\in ker\Psi_{G}^{_{\tau}G.P}$ for all $x\in CR(G,\, _{\tau}G.P) $. Thus $\tau(x)\in \mathrm{J}_{P}(\mathbb{F}_{p}(P))$ for all $x\in CR(G,\, _{\tau}G.P) $. The assertion is proved.
		Then there exists map $\lambda$ from $L^{n}\times U$ to $W_{n}$ such that 
		$$y=\sum\limits_{i_{1}+\cdots+i_{t}\geq1}\lambda(i_{1},\cdots,i_{t},1)e_{i_{1},\cdots,i_{t},1}.$$
		Next, verify uniqueness. Set $\sigma$ is a map $\lambda$ from $L^{n}\times U$ to $W_{n}$ that $$ \sum\limits_{i_{1}+\cdots+i_{t}\geq1}\sigma(i_{1},\cdots,i_{t},1)e_{i_{1},\cdots,i_{t},1}=0.$$
		Suppose $\sigma(i_{1},\cdots,i_{t},1)=\sum\limits_{x\in CR(G,\, _{\tau}G.P)}\d_{i_{1},\cdots,t_{t},x}x$. Then $$0=\sum\limits_{i_{1}+\cdots+i_{t}\geq1}\sigma(i_{1},\cdots,i_{t},1)e_{i_{1},\cdots,i_{t},1}=\sum\limits_{x\in CR(G,\, _{\tau}G.P)}x(\sum\limits_{i_{1}+\cdots+i_{t}\geq1}\d_{i_{1},\cdots,t_{t},x}e_{i_{1},\cdots,i_{t},1}).$$
		From $ \{(a_{1}-1)^{i_{1}}\cdots(a_{t}-1)^{i_{t}}|\sum\limits_{j=1}^{t}i_{j}\geq 1\}$ is a basis of $\mathrm{J}_{P}(\mathbb{F}_{p}(P))$, $\d_{i_{1},\cdots,t_{t},x}=0 $ for $\sum\limits_{j=1}^{n}i_{j}\geq1$ and $x\in CR(G,\, _{\tau}G.P)$. Thus $\mathrm{Im}\sigma=\{0\}$. The uniqueness is proved.
		
		The remaining part of lemma can be proved similarly.
		
		The proof is completed.  
	\end{pf}

	\begin{lem}\label{to}
		Set
		\begin{align*}
			_{n}\Psi_{ _{\tau}G.P}^{G}:\prod\limits^{n}\mathbb{F}_{p}(G)&\rightarrow \prod\limits^{n}\mathbb{F}_{p}(_{\tau}G.P)\\
			(g_{1},\cdots,g_{n})&\mapsto(\Psi_{ _{\tau}G.P}^{G}(g_{1}),\cdots,\Psi_{ _{\tau}G.P}^{G}(g_{n})).
		\end{align*}
		Set $x^{h}=x(\eta_{\tau,G}(h),\cdots,\eta_{\tau,G}(h))$ for all $x\in\prod\limits^{n} \mathbb{F}_{p}(G),h\in _{\tau}G.P$. $\prod\limits^{n} \mathbb{F}_{p}(G) $ is a right $_{\tau}G.P$ module. And $ \prod\limits^{n} \mathbb{F}_{p}(G)\cong \mathrm{Im} _{n}\Psi_{_{\tau}G.P}^{G}$ as right $ _{\tau}G.P$-module. Similarly, $ \prod\limits^{n} \mathbb{F}_{p}(G)\cong \mathrm{Im} _{n}\Psi_{_{\tau}G.P}^{G}$ as left $ _{\tau}G.P$-module.
	\end{lem}
	\begin{pf}
		Claim $ x(a_{1}-1)^{p-1}\cdots(a_{n}-1)^{p-1}=(a_{1}-1)^{p-1}\cdots(a_{n}-1)^{p-1}x$ for all $x\in \mathbb{F}_{p}(_{\tau}G.P)$.
		Notice $\sum\limits_{g\in P}g=(a_{1}-1)^{p-1}\cdots(a_{n}-1)^{p-1}$. From $xP=Px$ for all $x\in _{\tau}G.P$, $ x(a_{1}-1)^{p-1}\cdots(a_{n}-1)^{p-1}=(a_{1}-1)^{p-1}\cdots(a_{n}-1)^{p-1}x$ for all $x\in _{\tau}G.P$. The assertion is proved.
		
		Since $ \prod\limits^{n}\mathbb{F}_{p}(G)=\mathrm{Im}_{n}\Psi_{G}^{_{\tau }G.P}$ and $ker_{n}\Psi_{G}^{_{\tau }G.P}$ is a right $_{\tau}G.P$ -module,  $\prod\limits^{n}\mathbb{F}_{p}(G)\cong \prod\limits^{n}\mathbb{F}_{p}(_{\tau}G.P)/ker_{n}\Psi_{G}^{_{\tau }G.P}$ as $_{\tau}G.P $-module. Notice $\mathrm{Im} _{n}\Psi_{_{\tau}G.P}^{G}$ is right $_{\tau}G.P$-module generated by $e_{p-1,\cdots,p-1,1},\cdots,e_{p-1,\cdots,p-1,n}$. From $_{n}\Psi_{_{\tau}G.P}^{G} $ is injective, only
		verify $\prod\limits^{n}\mathbb{F}_{p}(_{\tau}G.P)/ker_{n}\Psi_{G}^{_{\tau }G.P}\cong \mathrm{Im} _{n}\Psi_{_{\tau}G.P}^{G}$ as right $_{\tau}G.P$-module. Set 
		\begin{align*}
			\psi:\prod\limits^{n}\mathbb{F}_{p}(_{\tau}G.P)/ker_{n}\Psi_{G}^{_{\tau }G.P}&\rightarrow \mathrm{Im} _{n}\Psi_{_{\tau}G.P}^{G}\\
			x+ker_{n}\Psi_{G}^{_{\tau }G.P}&\mapsto x\sum\limits_{i=1}^{n}e_{p-1,\cdots,p-1,i}.
		\end{align*}
		Since $ (a_{i}-1)^{p}=0$ for all $1\leq i\leq n$, $a_{i}$ acts on $\mathrm{Im} _{n}\Psi_{_{\tau}G.P}^{G}$ trivially for all $1\leq i\leq n$. Then $P$ acts on $ \mathrm{Im} _{n}\Psi_{_{\tau}G.P}^{G}$ trivially. Thus $\psi$ is well defined. Since $ _{n}\Psi_{G}^{_{\tau}G.P}$ is epimorphism, $\psi$ is epimorphism. And \begin{align*}\psi((x+ker_{n}\Psi_{G}^{_{\tau }G.P})^{y})&=\psi(x(y,\cdots,y)+ker_{n}\Psi_{G}^{_{\tau }G.P})\\&=x(y,\cdots,y)\sum\limits_{i=1}^{n}e_{p-1,\cdots,p-1,i}\\&=(x\sum\limits_{i=1}^{n}e_{p-1,\cdots,p-1,i})^{y}\end{align*}
		Then $\prod\limits^{n}\mathbb{F}_{p}(_{\tau}G.P)/ker_{n}\Psi_{G}^{_{\tau }G.P}\cong \mathrm{Im} _{n}\Psi_{_{\tau}G.P}^{G}$ as right $_{\tau}G.P$-module.
		
		The proof is completed.
	\end{pf}

	\begin{rem} Prove $\prod\limits^{n}\mathbb{F}_{p}(_{\tau}G.P)/ker_{n}\Psi_{G}^{_{\tau }G.P}\cong \mathrm{Im} _{n}\Psi_{_{\tau}G.P}^{G}$ as left $_{\tau}G.P$-module similarly.
		Set $_{n}\Psi_{_{\tau}G.P}^{_{\tau}G.P}= (_{n}\Psi_{_{\tau}G.P}^{G}) (_{n}\Psi_{G}^{_{\tau}G.P})$. Then $$_{n}\Psi_{_{\tau}G.P}^{_{\tau}G.P}(x)=x((a_{1}-1)^{p-1}\cdots(a_{t}-1)^{p-1},\cdots,(a_{1}-1)^{p-1}\cdots(a_{t}-1)^{p-1}) $$  for all $x\in \prod\limits^{n}\mathbb{F}_{p}(_{\tau}G.P)$.
	\end{rem}
	
	\begin{thm}\label{2.e}
		Let $ Q$ be a $G$-submodule of $\prod\limits^{n}\mathbb{F}_{p}(G)$. Denote $ H$ as image of $Q$ in $_{n}\Psi_{_{\tau}G.P}^{G}$. Then the image of $\mathcal{L}_{_{\tau}G.P}(H)$ in $_{n}\Psi_{G}^{_{\tau}G.P}$ is $\mathcal{L}_{G}(Q)$.
	\end{thm}
	\begin{pf}
		Without loss generality, assume $Q$ is a right $G$-submodule.
		
		Set $H_{1}=(_{n}\Psi_{_{\tau}G.P}^{_{\tau}G.P})^{-1}(H)$. Since $_{n}\Psi_{G}^{_{\tau}G.P}$ is surjective and $_{n}\Psi_{_{\tau}G.P}^{G}$ is injective, $_{n}\Psi_{G}^{_{\tau}G.P}(H_{1})=Q$. From $_{n}\Psi_{G}^{_{\tau}G.P}$ is surjective algebra homomorphism from $ \prod\limits^{n}\mathbb{F}_{p}(_{\tau}G.P)$ to $\prod\limits^{n}\mathbb{F}_{p}(G)$, only verify $$ \mathcal{L}_{_{\tau}G.P}(H)=\{x\in \prod\limits^{n}\mathbb{F}_{p}(_{\tau}G.P)|x\cdot y\in ker _{n}\Psi_{G}^{_{\tau}G.P}\; for\; all\; y\in H_{1}\}.$$
		Notice $(x\cdot y)(a_{1}-1)^{p-1}\cdots(a_{t}-1)^{p-1}=x\cdot _{n}\Psi_{_{\tau}G.P}^{_{\tau}G.P}(y)$ for all $x,y\in \prod\limits^{n}\mathbb{F}_{p}(_{\tau}G.P)$. Since $z(a_{1}-1)^{p-1}\cdots(a_{t}-1)^{p-1}=0$ if and only if $z\in ker\Psi_{G}^{_{\tau}G.P}$, $\{x\in \prod\limits^{n}\mathbb{F}_{p}(_{\tau}G.P)|x\cdot y\in ker\Psi_{G}^{_{\tau}G.P}\; for\; all\; y\in H_{1}\}=\{x\in \prod\limits^{n}\mathbb{F}_{p}(_{\tau}G.P)|x\cdot \Psi_{_{\tau}G.P}^{_{\tau}G.P}(y)=0\; for\; all\; y\in H_{1}\}. $ Notice $H_{1}=(_{n}\Psi_{_{\tau}G.P}^{_{\tau}G.P})^{-1}(H)$.  Thus the image of $\mathcal{L}_{_{\tau}G.P}(H)$ in $_{n}\Psi_{G}^{_{\tau}G.P}$ is $\mathcal{L}_{G}(Q)$.
		
		The proof is completed.
	\end{pf}
	
	\begin{defi}
		Let $x_{1},\cdots,x_{s}\in\prod\limits^{n}\mathbb{F}_{p}(G)$. Denote $$Ann_{L}(x_{1},\cdots,x_{s})=\{(y_{1},\cdots,y_{s})\in \prod\limits^{s}\mathbb{F}_{p}(G)|\sum\limits_{i=1}^{s}(\overbrace{y_{i},\cdots,y_{i}}^{n})x_{i}=0\}$$
		and $$Ann_{R}(x_{1},\cdots,x_{s})=\{(y_{1},\cdots,y_{s})\in \prod\limits^{s}\mathbb{F}_{p}(G)|\sum\limits_{i=1}^{s}x_{i}(\overbrace{y_{i},\cdots,y_{i}}^{n})=0\}$$
	\end{defi}

	Set $I_{n,m}$ is left $_{\tau}G.P$-module generated by $\{e_{i_{1},\cdots,i_{t},l}|\sum\limits_{j=1}^{t}i_{j}\geq m\}$.  And $I_{n,1}=ker _{n}\Psi_{G}^{_{\tau}G.P},I_{n,0}=\prod\limits^{n}\mathbb{F}_{p}(_{\tau}G.P)$.
	
	\begin{lem}\label{tp}
		$I_{n,m}$ is also right $_{\tau}G.P$-module generated by $\{e_{i_{1},\cdots,i_{t},l}|\sum\limits_{j=1}^{t}i_{j}\geq m\}$ for all $m\in\mathbb{N}_{+}$. 	$I_{n,m_{1}}I_{n,m_{2}}=I_{n,m_{1}+m_{2}}$.
	\end{lem}
	\begin{pf}Set $A_{m}=\{e_{i_{1},\cdots,i_{t},l}|\sum\limits_{j=1}^{t}i_{j}\geq m\}$.
		From Lemma \ref{xo}, only verify $e_{i_{1},\cdots,i_{t},l}x\in I_{n,m}$ for all $e_{i_{1},\cdots,i_{t},l}\in A_{m}, x\in W_{n}$. Since $ I_{n,1}=ker _{n}\Psi_{G}^{_{\tau}G.P}$, $I_{n,1}$ is left $_{\tau}G.P$-module generated by $A_{1}$. By induction on $m$.   When $m\geq 2$, there exists $e_{j_{1},\cdots,j_{t},l}\in A_{1}\backslash A_{2},e_{r_{1},\cdots,r_{t},l}\in A_{m-1}$ such that $e_{r_{1},\cdots,r_{t},l}e_{j_{1},\cdots,j_{t},l}= e_{i_{1},\cdots,i_{t},l}$. Since $I_{n,1}$ is left $_{\tau}G.P$-module generated by $A_{1}$, there exists map $\lambda$ from $ L^{t}$ to $ W_{n}$ such that $$e_{j_{1},\cdots,j_{t},l}x=\sum\limits_{q_{1}+\cdots+q_{t}\geq 1}\lambda(q_{1},\cdots,q_{t})e_{q_{1},\cdots,q_{t},l}.$$
		Since $e_{r_{1},\cdots,r_{t},l}\lambda(q_{1},\cdots,q_{t})\in I_{n,m-1}$ and $A_{m_{1}}A_{m_{2}}=A_{m_{1}+m_{2}}$, $e_{i_{1},\cdots,i_{t},l}x\in I_{n,m}$ for all $e_{i_{1},\cdots,i_{t},l}\in A_{m}, x\in W_{n}$. Thus  $I_{n,m}$ is right $_{\tau}G.P$-module. Set $ T$ is right $_{\tau}G.P$- module generated by $A_{m}$. Then $T\leq I_{n,m}$. It's similarly to prove $ T$ is a left $_{\tau}G.P$-module. Then $T= I_{n,m}$. Thus $I_{n,m}$ is also right $_{\tau}G.P$-module generated by $\{e_{i_{1},\cdots,i_{t},l}|\sum\limits_{j=1}^{t}i_{j}\geq m\}$ for all $m\in\mathbb{N}_{+}$. 
		
		From $I_{n,m}$ is also right $_{\tau}G.P$-module,  $I_{n,m_{1}}I_{n,1}=I_{n,m_{1}+1}$ for all $m_{1}\in \mathbb{N}_{+}$. Then $I_{n,m_{1}}I_{n,m_{2}}=I_{n,m_{1}+m_{2}}$.
		
		The proof is completed.
	\end{pf}

	\begin{lem}\label{dd}
		$\eta_{\tau,G}$ acts on $I_{n,i}/I_{n,i+1}$ trivially. Then $I_{n,1}/I_{n,2}$ is isomorphic to direct sum of $nt$ copies of $\mathbb{F}_{p}(G)$ as left(right) $G$-module. And $ker_{n}\Psi_{G}^{_{\tau}G.P}$ has $_{\tau}G.P $ minimal generator set of order $nt$.
	\end{lem}
	\begin{pf}Only discuss the aspect of $I_{n,i}$ as left $_{\tau}G.P$-module. 
		From Lemma \ref{tp}, $I_{n,1}I_{n,i}=I_{n,i}I_{n,1}= I_{n,i+1}$.
		Thus $\eta_{\tau,G}$ acts on $I_{n,i}/I_{n,i+1}$ trivially. And $I_{n,i}/I_{n,i+1}$ can be recorded as $G$-module. From Lemma \ref{xo}, $I_{n,1}/I_{n,2}$ is isomorphic to direct sum of $nt$ copies of $\mathbb{F}_{p}(G)$ as left(right) $G$-module. Since $I_{n,2}=I_{n,1}I_{n,1}$, $I_{n,2}\leq\mathrm{J}_{_{\tau}G.P}(I_{n,1})$. Thus $d_{_{\tau}G.P}(I_{n,1})=d_{G}(I_{n,1}/I_{n,2})$. And $ker_{n}\Psi_{G}^{_{\tau}G.P}$ has $_{\tau}G.P $ minimal generator set of order $nt$.
		
		The proof is completed.
	\end{pf}
	
	\begin{rem}\label{xp}
		When $t=2$, $I_{n,i}/I_{n,i+1}$ is isomorphic to direct sum of $ i+1$ copies of $ \mathbb{F}_{p}(G)$ as right(left) $ G$-module where $ i\leq p-1$. And $I_{n,i}/I_{n,i+1}$ is isomorphic to direct sum of $ 2p-1-i$ copies of $ \mathbb{F}_{p}(G)$ as right(left) $ G$-module where $ p\leq i\leq 2p-2$.
	\end{rem}
	
	\begin{cor}\label{yy}
		Let $Q$ be a $n-G$-module. Then $ \mathrm{H}^{1}(_{\tau}G.P,Q)\leq \mathrm{H}^{1}(G,Q)\oplus\mathbb{F}_{p}^{nt}$.
	\end{cor}
	\begin{pf}
		Without loss of generality, suppose $Q$ is right $G$-module. Set $Q_{1}$ is right $G$-submodule of $\prod\limits^{n}\mathbb{F}_{p}(G)$ that $Q\cong Q_{1}$. Denote $Q_{2}=\, _{n}\Psi^{G}_{_{\tau}G.P}(Q_{1})$. Suppose $\mathrm{H}^{1}(G,Q)\cong\mathbb{F}_{p}^{m}$. Then $\mathcal{L}_{G}(Q_{1})$ has $G$ minimal generator set of order $m$. From Theorem \ref{2.e} and Lemma \ref{dd}, $\mathcal{L}_{_{\tau}G.P}(Q_{2})$ has $_{\tau}G.P$ minimal generator set of order at most $m+nt$. From Lemma \ref{ww},
		$ \mathrm{H}^{1}(_{\tau}G.P,Q)\leq \mathrm{H}^{1}(G,Q)\oplus\mathbb{F}_{p}^{nt}$.
		
		The proof is completed.
	\end{pf}

	\begin{lem}\label{tu}
		Let $Q$ be a left $_{\tau}G.P $-submodule of $\prod\limits^{n}\mathbb{F}_{p}(_{\tau}G.P)$ that $d_{G}(_{n}\Psi_{G}^{_{\tau}G.P }(Q))=d_{_{\tau}G.P}(Q)\leq n$. Denote $\{x_{1},\cdots,x_{m}\}$ as a $_{\tau}G.P$ minimal generator set of $Q$. Set $ g_{i}=(_{n}\Psi_{G}^{_{\tau}G.P })(x_{i})$. Set $$D=\{(b_{1},\cdots,b_{m})\in \prod\limits^{m}\mathbb{F}_{p}(_{\tau}G.P)|\sum\limits_{i=1}\limits^{m}(\overbrace{b_{i},\cdots,b_{i}}^{n})x_{i}\in ker_{n}\Psi_{G}^{_{\tau}G.P}\}$$ and
		\begin{align*}
			\xi:D&\rightarrow ker _{n}\Psi_{G}^{_{\tau}G.P }/I_{n,2}\\
			(b_{1},\cdots,b_{m})&\mapsto \sum\limits_{i=1}\limits^{m}(b_{i},\cdots,b_{i})x_{i}+I_{n,2}.
		\end{align*}
		Then $ |coker\xi|\geq \frac{|\mathbb{F}_{p}(G)|^{nt-m}}{|_{n}\Psi_{G}^{_{\tau}G.P }(Q)|^{t-1}}$.
	\end{lem}
	\begin{pf}
		Denote $D_{1}=Ann_{L}(g_{1},\cdots,g_{m})$. Set \begin{align*}
			X:\prod\limits^{m}\mathbb{F}_{p}(_{\tau}G.P)&\rightarrow \prod\limits^{m}\mathbb{F}_{p}(G)\\
			(u_{1},\cdots,u_{m})&\mapsto (\Psi_{G}^{_{\tau}G.P }(u_{1}),\cdots,\Psi_{G}^{_{\tau}G.P }(u_{m})).
		\end{align*}
		Notice $ X(D)=D_{1}$ and $kerX= ker _{m}\Psi_{G}^{_{\tau}G.P } $. Since $|_{n}\Psi_{G}^{_{\tau}G.P }(Q)||D_{1}|=|\mathbb{F}_{p}(G)|^{m} $, $|D/kerX|=|D_{1}|=\frac{ |\mathbb{F}_{p}(G)|^{m}}{|_{n}\Psi_{G}^{_{\tau}G.P }(Q)|}$.
		Denote $ e(x,l)$ as element of $\prod\limits^{m}\mathbb{F}_{p}(_{\tau}G.P)$ that $l.th$ component is $x$ and others components are zero where $ x\in CR(G,\,_{\tau}G.P)$. Set $E=\{e(x,l)|x\in CR(G,\,_{\tau}G.P),1\leq l\leq m\}$ and $E_{1}$ is $\mathbb{F}_{p}$ linear space over basis $E$. 
		
		Set $B_{i}=\{((a_{i}-1)u_{1},\cdots,(a_{i}-1)u_{m})|(u_{1},\cdots,u_{m})\in E_{1}\}$ and 
		\begin{align*}
			\lambda_{i}: E_{1}&\rightarrow B_{i}\\
			(u_{1},\cdots,u_{m})&\mapsto ((a_{i}-1)u_{1},\cdots,(a_{i}-1)u_{m}).
		\end{align*}
		Notice $ \lambda_{i}$ is bijective.
		For $(u_{1},\cdots,u_{m})\in B_{i} $, $ (u_{1},\cdots,u_{m})\in ker\xi $ if and only if $ ((\lambda_{i})^{-1}(u_{1}),\cdots,(\lambda_{i})^{-1}(u_{m}))\in E_{1}\cap D$. From Lemma \ref{xo}, 
		$$I_{m,1}=B_{1}\oplus B_{2}\oplus\cdots\oplus B_{t}\oplus I_{m,2}.$$
		Since $ I_{m,2}$ is a two-sided $_{\tau}G.P $-submodule of $\prod\limits^{m}\mathbb{F}_{p}(_{\tau}G.P)$, $I_{m,2}\leq ker\xi$.
		Notice $|E_{1}\cap D|=| D/kerX|$. Thus $|\xi(kerX)|=(\frac{ |E_{1}|}{|D/kerX|})^{t}=\frac{|\mathbb{F}_{p}(G)|^{mt}}{|D/kerX|^{t}}$. And $$ |\mathrm{Im}\xi|\leq |\xi(kerX)||D/kerX|=\frac{|\mathbb{F}_{p}(G)|^{mt}}{|D/kerX|^{t-1}}.$$
		From $|D/kerX|=\frac{ |\mathbb{F}_{p}(G)|^{m}}{|_{n}\Psi_{G}^{_{\tau}G.P }(Q)|}$, $$|\mathrm{Im}\xi|\leq |\mathbb{F}_{p}(G)|^{m}|_{n}\Psi_{G}^{_{\tau}G.P }(Q)|^{t-1}.$$
		From Lemma \ref{dd},  $ |coker\xi|\geq \frac{|\mathbb{F}_{p}(G)|^{nt-m}}{|_{n}\Psi_{G}^{_{\tau}G.P }(Q)|^{t-1}}$.
		
		The proof is completed.
	\end{pf}
	
	\begin{thm}\label{jj}
		Let $ Q$ be a right(left) $_{\tau}G.\mathrm{C}_{p}$-submodule of $\prod\limits^{n}\mathbb{F}_{p}(_{\tau}G.\mathrm{C}_{p})$ that  $ Q\leq \mathrm{Im}_{n}\Psi^{G}_{_{\tau}G.\mathrm{C}_{p}}$ and $ \mathrm{H}^{1}(G,Q)\cong\mathbb{F}_{p}^{m}$ where $m<n$. Then $\mathcal{L}_{_{\tau}G.\mathrm{C}_{p}}(Q)(\mathcal{R}_{_{\tau}G.\mathrm{C}_{p}}(Q))$ has a $_{\tau}G.\mathrm{C}_{p} $ minimal generator set of order $s$ where $s\geq n$.
	\end{thm}
	\begin{pf} Without loss of generality, suppose $Q$ is right ${\tau}G.\mathrm{C}_{p}$-module.
		Denote $ Q_{1}= (_{n}\Psi^{G}_{_{\tau}G.\mathrm{C}_{p}})^{-1}(Q)$.	Since $ \mathrm{H}^{1}(G,Q)\cong \mathbb{F}_{p}^{m}$, $\mathcal{L}_{G}(Q_{1})$ has a $G$ minimal generator set of order $m$ as left $\mathbb{F}_{p}(G)$-module. Set $g_{1},\cdots,g_{m}\in \prod\limits^{n}\mathbb{F}_{p}(G)$ that $ \mathcal{L}_{G}( Q_{1})$ is left $G$-submodule of $\prod\limits^{n}\mathbb{F}_{p}(G)$ generated by $g_{1},\cdots,g_{m}$. Denote $x_{1},\cdots,x_{m}\in \prod\limits^{n}\mathbb{F}_{p}(_{\tau}G.\mathrm{C}_{p})$ that $_{n}\Psi_{G}^{_{\tau}G.\mathrm{C}_{p}}(x_{i})=g_{i}$ for all $1\leq i\leq m$. Then $x_{i}\in \mathcal{L}_{_{\tau}G.\mathrm{C}_{p}}(Q)$ for all $1\leq i\leq n$. Denote $T$ as left $_{\tau}G.\mathrm{C}_{p}$ module generated by $ x_{1},\cdots,x_{m}$.
		
		From Theorem \ref{2.e}, $\mathcal{L}_{_{\tau}G.\mathrm{C}_{p}}(Q)$ is left $_{\tau}G.\mathrm{C}_{p} $-module generated by $x_{1},\cdots,x_{m},I_{n,1}$. Assume $I_{n,1}/(T\cap I_{n,1})$ has $G$ minimal generator set of order at least $n-m$. From Lemma \ref{cc}, $\mathcal{L}_{_{\tau}G.\mathrm{C}_{p}}(Q) $ has $ _{\tau}G.\mathrm{C}_{p}$ minimal generator set of order at least $n$. Next, prove the assumption is true. 
		
		From Lemma \ref{tu}, $|ker_{n}\Psi^{_{\tau}G.\mathrm{C}_{p}}_{G}/(T+I_{n,2})|\geq|\mathbb{F}_{p}(G)|^{n-m}$. Since a $\mathbb{F}_{p}(G)$- module generated by $s$ elements has order at most $|\mathbb{F}_{p}(G)|^{s}$, 
		$ker_{n}\Psi^{_{\tau}G.\mathrm{C}_{p}}_{G}/(T+I_{n,2})$ has $G$ minimal genrator set of order at least $n-m$. Notice $ker_{n}\Psi^{_{\tau}G.\mathrm{C}_{p}}_{G}/(T+I_{n,2}) $ is quotient module of $ker_{n}\Psi^{_{\tau}G.\mathrm{C}_{p}}_{G}/(T\cap ker_{n}\Psi^{_{\tau}G.\mathrm{C}_{p}}_{G})$. And 
		$ker_{n}\Psi^{_{\tau}G.\mathrm{C}_{p}}_{G}/(T\cap ker_{n}\Psi^{_{\tau}G.\mathrm{C}_{p}}_{G})$ has $G$ minimal generator set of order at least $n-m$.
		
		The proof is completed.
	\end{pf}
	
	{\noindent\bf Proof of Theorem \ref{gg}}.
	
	\begin{pf}
		From Lemma \ref{ut}, $Q$ can be embedded in $\prod\limits^{n}\mathbb{F}_{p}(G)$ as $G$-module. Set $Q_{1}$ is $G$-submodule of $\prod\limits^{n}\mathbb{F}_{p}(G)$ that $Q_{1}\cong Q$ as $G$-module. Denote $ Q_{2}= _{n}\Psi_{_{\tau}G.\mathrm{C}_{p}}^{G}(Q_{1})$.
		
		From Theorem \ref{jj} and Lemma \ref{ww}, $\mathrm{H}^{1}(_{\tau}G.\mathrm{C}_{p},Q_{2})\geq\mathbb{F}_{p}^{n}$. Now, $$\mathrm{H}^{1}(_{\tau}G.\mathrm{C}_{p},Q)\geq \mathrm{H}^{1}(G,Q)\oplus\mathbb{F}_{p}^{n-m}.$$
	\end{pf}
	
	\begin{cor}\label{8.0}
		Denote $A$ as a $n-G$-module. Let $\tau\in Z^{2}(G,\mathrm{C}_{p})$. If\; $ \mathrm{H}^{1}(_{\tau}G.\mathrm{C}_{p},A)\cong\mathrm{H}^{1}(G,A)$, then $\mathrm{H}^{1}(G,A)\cong\mathbb{F}_{p}^{n}$.
	\end{cor}

	\begin{thm}\label{aa}
		Let $ Q$ be a right(left) ${\tau}G.P$-submodule of $\prod\limits^{n}\mathbb{F}_{p}(_{\tau}G.P)$ that  $ Q\leq \mathrm{Im}_{n}\Psi^{G}_{_{\tau}G.P}$ and $ \mathrm{H}^{1}(G,Q)\cong\mathbb{F}_{p}^{m}$ where $ m<n$. Then 
		
		{\rm(1)} $\mathrm{H}^{1}(_{\tau}G.P,Q)\geq\mathrm{H}^{1}(G,Q)\oplus\mathbb{F}_{p}^{tn-tm+1} $ when $m\geq 1$,
		
		{\rm(2)} $\mathrm{H}^{1}(_{\tau}G.P,Q)\cong\mathrm{H}^{1}(G,Q)\oplus\mathbb{F}_{p}^{tn} $ when $m=0$.
	\end{thm}
	\begin{pf}
		Without loss of generality, assume $ Q$ is a right $_{\tau}G.P$-submodule of $\prod\limits^{n}\mathbb{F}_{p}(_{\tau}G.P)$. Denote $Q_{1}=(_{n}\Psi^{G}_{_{\tau}G.P})^{-1}(Q)$.
		
		When $m=0$, $Q_{1}=\prod\limits^{n}\mathbb{F}_{p}(G)$. Then $ \mathcal{L}_{G}(Q_{1})=0$. From Lemma \ref{dd} and $\mathcal{L}_{_{\tau}G.P}(Q)=ker_{n}\Psi^{G}_{_{\tau}G.P}$, $ \mathcal{L}_{_{\tau}G.P}(Q)$ has $_{\tau}G.P $ minimal generator set of order $2nt$.
		
		Next, suppose $ m\geq 1$.
		Since $ \mathrm{H}^{1}(G,Q)\cong \mathbb{F}_{p}^{m}$, $\mathcal{L}_{G}((Q_{1})$ has a minimal generator set of order $m$ as left $\mathbb{F}_{p}(G)$-module. Set $g_{1},\cdots,g_{m}\in \prod\limits^{n}\mathbb{F}_{p}(G)$ that $ \mathcal{L}_{G}(Q_{1})$ is left $G$-submodule of $\prod\limits^{n}\mathbb{F}_{p}(G)$ generated by $g_{1},\cdots,g_{m}$. Denote $x_{1},\cdots,x_{m}\in \prod\limits^{n}\mathbb{F}_{p}(_{\tau}G.P)$ that $_{n}\Psi^{_{\tau}G.P}_{G}(x_{i})=g_{i}$ for all $1\leq i\leq m$. Then $x_{i}\in \mathcal{L}_{_{\tau}G.P}(Q)$ for all $1\leq i\leq n$. Denote $T$ as left $_{\tau}G.P$-submodule generated by $ x_{1},\cdots,x_{m}$.
		
		From Lemma \ref{tu}, $$| ker_{n}\Psi^{_{\tau}G.P}_{G}/(T+I_{n,2})|\geq \frac{|\mathbb{F}_{p}(G)|^{tn-m}}{|\mathcal{L}_{G}(Q_{1})|^{t-1}}.$$
		Since $ d_{G}(\mathcal{L}_{G}(Q_{1}))= m$ where $m\geq 1$, $|\mathcal{L}_{G}(Q_{1})|<|\mathbb{F}_{p}(G)|^{m}$. And $$| ker_{n}\Psi^{_{\tau}G.P}_{G}/(T+I_{n,2})|> |\mathbb{F}_{p}(G)|^{tn-tm} .$$
		Thus $ker_{n}\Psi^{_{\tau}G.P}_{G}/(T+I_{n,2})$ has $G$ minimal generator set of order at least $tn-tm+1$. Since $ker_{n}\Psi^{_{\tau}G.P}_{G}/(T+I_{n,2})$ is a quotient module of  $ker_{n}\Psi^{_{\tau}G.P}_{G}/(ker_{n}\Psi^{_{\tau}G.P}_{G}\cap T)$, $ker_{n}\Psi^{_{\tau}G.P}_{G}/(ker_{n}\Psi^{_{\tau}G.P}_{G}\cap T)$ has $_{\tau}G.P$ minimal generator set of order at least $tn-tm+1$.
		
		From Lemma \ref{ww}, the proof is completed.
	\end{pf}
	
	The above  discussion is about that $\mathrm{H}^{1}(G,Q)<\mathbb{F}_{p}^{n}$ where $Q^{G}\cong\mathbb{F}_{p}^{n}$. Next, the disscussion is about exactly $n-G$ module.
	
	\begin{lem}\label{q}
		Let $ Q$ be a right $G$-submodule of $\prod\limits^{n}\mathbb{F}_{p}(G)$ that $Q^{G}\cong\mathbb{F}_{p}^{n}$. Then $\mathcal{L}_{G}(Q)\leq \mathrm{J}_{G}(\prod\limits^{n}\mathbb{F}_{p}(G))$. 
	\end{lem}
	\begin{pf}
		If $\mathcal{L}_{G}(Q)\nleqslant \mathrm{J}_{G}(\prod\limits^{n}\mathbb{F}_{p}(G))$, there exists $(a_{1},\cdots,a_{n})\in \prod\limits^{n}\mathbb{F}_{p}(G)\backslash \mathrm{J}_{G}(\prod\limits^{n}\mathbb{F}_{p}(G))$. And $a_{i}\in\mathbb{F}_{p}(G)\backslash \mathrm{J}_{G}(\mathbb{F}_{p}(G))$ for some $1\leq i\leq n$. Set $h\in \prod\limits^{n}\mathbb{F}_{p}(G)$ that $ith$ component is $\sum\limits_{g\in G}g$ and others  components are zero. Then $h\notin\mathcal{R}_{G}(\mathcal{L}_{G}(Q))=Q$. It's contradictory to $h\in Q^{G}$.
		
		The proof is completed. 
	\end{pf}

	\begin{lem}\label{ii}
		Let $Q$ be a right $G$-submodule of $\prod\limits^{n}\mathbb{F}_{p}(G)$ that $\mathrm{H}^{1}(G,Q)\cong\mathrm{H}^{1}(_{\tau}G.\mathrm{C}_{p},Q)$ for $\tau\in Z^{2}(G,\mathrm{C}_{p})\backslash B^{2}(G,\mathrm{C}_{p}) $. Denote $Q_{1}=\,_{n}\Psi_{_{\tau}G.\mathrm{C}_{p}}^{G}(Q)$.  Set $x_{1},\cdots,x_{n}\in \prod\limits^{n}\mathbb{F}_{p}(_{\tau}G.\mathrm{C}_{p})$ that $ \mathcal{L}_{_{\tau}G.\mathrm{C}_{p}}(Q_{1})$ is left $_{\tau}G.\mathrm{C}_{p} $-module generated by $x_{1},\cdots,x_{n}$. Denote $D_{_{\tau}G.\mathrm{C}_{p}}=\{(y_{1},\cdots,y_{n})\in\prod\limits^{n}\mathbb{F}_{p}(_{\tau}G.\mathrm{C}_{p})|\sum\limits_{i=1}^{n}(\overbrace{y_{i},\cdots,y_{i}}^{n})x_{i}\in ker_{n}\Psi^{_{\tau}G.\mathrm{C}_{p}}_{G}\}$. Then $ D_{_{\tau}G.\mathrm{C}_{p}}$ is a left $_{\tau}G.\mathrm{C}_{p}$ module that $d_{_{\tau}G.\mathrm{C}_{p}}(D_{_{\tau}G.\mathrm{C}_{p}})=n$. And $ _{n}\Psi_{_{\tau}G.\mathrm{C}_{p}}^{_{\tau}G.\mathrm{C}_{p}}(D_{_{\tau}G.\mathrm{C}_{p}})=Ann_{L}(x_{1},\cdots,x_{n})$. Furthermore, $d_{_{\tau}G.\mathrm{C}_{p}}(Ann_{L}(x_{1},\cdots,x_{n}))=n$.
	\end{lem}
	\begin{pf} 
		Since $\mathrm{H}^{1}(G,Q)\cong\mathrm{H}^{1}(_{\tau}G.\mathrm{C}_{p},Q)$, $x_{1},\cdots,x_{n}\in \prod\limits^{n}\mathbb{F}_{p}(_{\tau}G.\mathrm{C}_{p})\backslash ker_{n}\Psi_{G}^{_{\tau}G.\mathrm{C}_{p}}$. From Lemma \ref{q}, $x_{1},\cdots,x_{n}\in \mathrm{J}_{G}(\prod\limits^{n}\mathbb{F}_{p}(G))$. Denote $h_{i}=\,_{n}\Psi_{G}^{_{\tau}G.\mathrm{C}_{p}}(x_{i})$. Then $\mathcal{L}_{G}(Q)$ is left $G$ module generated by $ h_{1},\cdots,h_{n}$. 
		
		Since $ |Ann_{L}(x_{1},\cdots,x_{n})||\mathcal{L}_{_{\tau}G.\mathrm{C}_{p}}|=|\prod\limits^{n}\mathbb{F}_{p}(_{\tau}G.\mathrm{C}_{p})|$, $$ |Ann_{L}(x_{1},\cdots,x_{n})|=|Ann_{L}(h_{1},\cdots,h_{n})|.$$ Then $Ann_{L}(x_{1},\cdots,x_{n})=\,_{n}\Psi^{G}_{_{\tau}G.\mathrm{C}_{p}}(Ann_{L}(h_{1},\cdots,h_{n}))$. Thus $ _{n}\Psi_{_{\tau}G.\mathrm{C}_{p}}^{_{\tau}G.\mathrm{C}_{p}}(D_{_{\tau}G.\mathrm{C}_{p}})=Ann_{L}(x_{1},\cdots,x_{n})$. And $Ann_{L}(x_{1},\cdots,x_{n})\leq \mathrm{J}_{_{\tau}G.\mathrm{C}_{p}}(\mathrm{Im}_{n}\Psi^{G}_{_{\tau}}G.\mathrm{C}_{p}) $.
		
		Notice $\sum\limits_{i=1}^{n}(\overbrace{y_{i},\cdots,y_{i}}^{n})x_{i}\in\mathrm{J}_{_{\tau}G.\mathrm{C}_{p}}(ker_{n}\Psi_{G}^{_{\tau}G.\mathrm{C}_{p}})$ for all $(y_{1},\cdots,y_{n})\in \mathrm{J}_{_{\tau}G.\mathrm{C}_{p}}(D_{_{\tau}G.\mathrm{C}_{p}})$ where $ ker_{n}\Psi_{G}^{_{\tau}G.\mathrm{C}_{p}}$ is recorded as left $_{\tau}G.\mathrm{C}_{p}$-module. Set 
		\begin{align*}
			\psi:D_{_{\tau}G.\mathrm{C}_{p}}/\mathrm{J}_{_{\tau}G.\mathrm{C}_{p}}(D_{_{\tau}G.\mathrm{C}_{p}})&\rightarrow ker_{n}\Psi_{G}^{_{\tau}G.\mathrm{C}_{p}}/ \mathrm{J}_{_{\tau}G.\mathrm{C}_{p}}(ker_{n}\Psi_{G}^{_{\tau}G.\mathrm{C}_{p}})\\
			(y_{1},\cdots,y_{n})+\mathrm{J}_{_{\tau}G.\mathrm{C}_{p}}(D_{_{\tau}G.\mathrm{C}_{p}})&\mapsto \sum\limits_{i=1}^{n}(\overbrace{y_{i},\cdots,y_{i}}^{n})x_{i}+\mathrm{J}_{_{\tau}G.\mathrm{C}_{p}}(ker_{n}\Psi_{G}^{_{\tau}G.\mathrm{C}_{p}}).
		\end{align*} 
		Since $ \mathcal{L}_{_{\tau}G.\mathrm{C}_{p}}(Q_{1})$ is left $_{\tau}G.\mathrm{C}_{p} $-module generated by $x_{1},\cdots,x_{n}$, $\psi$ is surjective. Then $ d_{_{\tau}G.\mathrm{C}_{p}}(D_{_{\tau}G.\mathrm{C}_{p}})\geq n$. 
		
		Set $a\in\,_{\tau}G.\mathrm{C}_{p}\backslash\{1\}$ that $a-1\in ker\Psi_{G}^{_{\tau}G.\mathrm{C}_{p}}$. Denote $e_{i}$ as element of $\prod\limits^{n}\mathbb{F}_{p}(_{\tau}G.\mathrm{C}_{p})$ that $i.th$ component is $a-1$ and others components are zero. Then there exists $z_{i}=(z_{1,i},\cdots,z_{n,i})$ such that $\sum\limits_{j=1}^{n}(\overbrace{z_{j,i},\cdots,z_{j,i}}^{n})x_{i}=e_{i}$. 
		
		Denote $H$ as left $_{\tau}G.\mathrm{C}_{p}$-module generated by $ z_{1},\cdots,z_{n}$.	If $d_{_{\tau}G.\mathrm{C}_{p}}(D_{_{\tau}G.\mathrm{C}_{p}})>n$, there exists $ u=(u_{1},\cdots,u_{n})\in D_{_{\tau}G.\mathrm{C}_{p}}\backslash (H+\mathrm{J}_{_{\tau}G.\mathrm{C}_{p}}(D_{_{\tau}G.\mathrm{C}_{p}}))$. 
		Since $ker_{n}\Psi_{G}^{_{\tau}G.\mathrm{C}_{p}} $ is left $_{\tau}G.\mathrm{C}_{p}$-module generated by $e_{1},\cdots,e_{n}$, there exists $(r_{1},\cdots,r_{n}) \in \prod\limits^{n}\mathbb{F}_{p}(_{\tau}G.\mathrm{C}_{p})$ such that $ u+\sum\limits_{j=1}^{n}(\overbrace{r_{j},\cdots,r_{j}}^{n})z_{j}\in Ann_{L}(x_{1},\cdots,x_{n})$.  Notice $u=(u_{1},\cdots,u_{n})\in D_{_{\tau}G.\mathrm{C}_{p}}\backslash (H+\mathrm{J}_{_{\tau}G.\mathrm{C}_{p}}(D_{_{\tau}G.\mathrm{C}_{p}}))$ . Then $u\notin \mathrm{J}_{_{\tau}G.\mathrm{C}_{p}}(ker_{n}\Psi_{G}^{_{\tau}G.\mathrm{C}_{p}})$. It's contradictory to $$Ann_{L}(x_{1},\cdots,x_{n})\leq \mathrm{J}_{_{\tau}G.\mathrm{C}_{p}}(\mathrm{Im}_{n}\Psi_{_{\tau}G.\mathrm{C}_{p}}^{G})\leq  \mathrm{J}_{_{\tau}G.\mathrm{C}_{p}}(ker_{n}\Psi_{G}^{_{\tau}G.\mathrm{C}_{p}}).$$ Then $H=D_{_{\tau}G.\mathrm{C}_{p}}.$
		
		From  $x_{1},\cdots,x_{n}\in \mathrm{J}_{G}(\prod\limits^{n}\mathbb{F}_{p}(G))$, $\sum\limits_{i=1}^{n}(\overbrace{y_{i},\cdots,y_{i}}^{n})x_{i}\in \mathrm{J}_{_{\tau}G.\mathrm{C}_{p}}(ker_{n}\Psi_{G}^{_{\tau}G.\mathrm{C}_{p}})$ for all $(y_{1},\cdots,y_{n})\in ker_{n}\Psi_{G}^{_{\tau}G.\mathrm{C}_{p}}$. Then $ \psi(ker_{n}\Psi_{G}^{_{\tau}G.\mathrm{C}_{p}})=0 $. Notice $Ann_{L}(x_{1},\cdots,x_{n})\cong D_{_{\tau}G.\mathrm{C}_{p}}/ker_{n}\Psi_{G}^{_{\tau}G.\mathrm{C}_{p}}$ as left $_{\tau}G.\mathrm{C}_{p} $ module. Then $Ann_{L}(x_{1},\cdots,x_{n}) $ has  $_{\tau}G.\mathrm{C}_{p} $ minimal generator set  $$_{n}\Psi_{_{\tau}G.\mathrm{C}_{p}}^{_{\tau}G.\mathrm{C}_{p}}(z_{1}),\cdots,\,_{n}\Psi_{_{\tau}G.\mathrm{C}_{p}}^{_{\tau}G.\mathrm{C}_{p}}(z_{n})$$ as left $_{\tau}G.\mathrm{C}_{p} $-module.
		
		The proof is completed.
	\end{pf}

	\begin{rem}\label{rt}
		When $Q$ is left module, there are results similar to the above two lemmas.
	\end{rem}

	\begin{thm}\label{ggg}
		Let $Q$ be a right(left) exactly $n-G$-submodule of $\prod\limits^{n}\mathbb{F}_{p}(G)$. Denote $ N$ as $\mathbb{F}_{p}(G)$-module that $|N|=p^{2}$. $N_{1}$ is a $G$-submodule of $ N$ that $|N_{1}|=p$. If exists $\tau\in Z^{2}(G,N)$ such that $ \mathrm{H}^{1}(_{\tau}G.N/N_{1},Q)=\mathrm{H}^{1}(G,Q)$, then $\mathrm{H}^{1}(_{\tau}G.N,Q)\cong \mathrm{H}^{1}(_{\tau}G.N/N_{1},Q)\oplus\mathrm{F}_{p}^{n}$.
	\end{thm}
	\begin{pf}Without loss of generality, suppose $Q$ is a right exactly $n-G$-submodule of $\prod\limits^{n}\mathbb{F}_{p}(G)$.
		Suppose $ N=\langle a_{1}\rangle\times \langle a_{2}\rangle$ and $N_{1}=\langle a_{2}\rangle$. Denote $H=\,_{\tau}G.N/N_{1}$. And $\overline{a_{1}}$ is the image of $a_{1}$ in natural group epimorphism from $_{\tau}G.N$ to $H$. $\overline{1}$ is identity element of $H$. Set $Q_{1}=\,_{n}\Psi_{G}^{H}(Q)$ and $Q_{2}=\,_{n}\Psi_{G}^{_{\tau}G.N}(Q)$.
		
		From $ \mathrm{H}^{1}(_{\tau}G.N/N_{1},Q)=\mathrm{H}^{1}(G,Q)$, $\mathcal{L}_{H}(Q_{1})$ is left $H$-module that $d_{H}(\mathcal{L}_{H}(Q_{1}))=n$. Suppose  $\mathcal{L}_{H}(Q_{1}) $ is left $H$-module generated by $x_{1},\cdots,x_{n}\in \mathrm{J}_{H}(\prod\limits^{n}\mathbb{F}_{p}(H))\backslash  ker_{n}\Psi_{G}^{H}$. Denote $\overline{e}_{i}=((\overline{a_{1}}-\overline{1})^{i},\cdots,(\overline{a_{1}}-\overline{1})^{i})$ and $e_{i}=((a_{1}-1)^{i},\cdots,(a_{1}-1)^{i})$. From Lemma \ref{ii}, there exist $\overline{z}_{1},\cdots,\overline{z}_{n}\in \mathrm{J}_{H}(\prod\limits^{n}\mathbb{F}_{p}(H))\backslash  ker_{n}\Psi_{G}^{H}$ such that $ Ann_{L}(x_{1},\cdots,x_{n})$ is left $_{\tau}G.\mathrm{C}_{p}$-module generated by $\overline{z}_{1}\overline{e}_{p-1},\cdots,\overline{z}_{n}\overline{e}_{p-1}$.
		
		Set $y_{1},\cdots,y_{n}\in \prod\limits^{n}\mathbb{F}_{p}(_{\tau}G.N)$ that $_{n}\Psi_{H}^{_{\tau}G.N}(y_{i})=x_{i}$. And $z_{1},\cdots,z_{n}\in \prod\limits^{n}\mathbb{F}_{p}(_{\tau}G.N)$ that $_{n}\Psi_{H}^{_{\tau}G.N}(z_{i})=\overline{z}_{i}$. Denote $A$ as left $_{\tau}G.N$ -module generated by $ y_{1},\cdots,y_{n}$. Assume $A\cap ker_{n}\Psi_{H}^{_{\tau}G.N}\leq \mathrm{J}_{_{\tau}G.N}(_{n}\Psi_{H}^{_{\tau}G.N})$. Then $\mathcal{L}_{_{\tau}G.N}(Q_{2})$ has $_{\tau}G.N$ minimal generator set of order $2n$. Next, prove the assumption.
		
		Set $D_{_{\tau}G.N}=\{(u_{1},\cdots,u_{n})|\sum\limits_{i=1}^{n}(u_{i},\cdots,u_{i})y_{i}\in ker_{n}\Psi_{H}^{_{\tau}G.N}\}$. Then $D_{_{\tau}G.N}$ is left $_{\tau}G.N$-module generated by $z_{1}e_{p-1},\cdots,z_{n}e_{p-1},ker_{n}\Psi_{H}^{_{\tau}G.N}$. Set
		\begin{align*}
			\phi:D_{_{\tau}G.N}&\rightarrow ker_{n}\Psi_{H}^{_{\tau}G.N}\\
			(u_{1},\cdots,u_{n})&\mapsto \sum\limits_{i=1}^{n}(u_{i},\cdots,u_{i})y_{i}.
		\end{align*} When $\mathrm{Im}\phi\leq \mathrm{J}_{_{\tau}G.N}(ker_{n}\Psi_{H}^{_{\tau}G.N})$, $A\cap ker_{n}\Psi_{H}^{_{\tau}G.N}\leq \mathrm{J}_{_{\tau}G.N}(ker_{n}\Psi_{H}^{_{\tau}G.N})$.
		From $x_{1},\cdots,x_{n}\in \mathrm{J}_{H}(\prod\limits^{n}\mathbb{F}_{p}(H))$, $y_{1},\cdots,y_{n}\in \mathrm{J}_{_{\tau}G.N}(\prod\limits^{n}\mathbb{F}_{p}(_{\tau}G.N))$. Then $$\phi(ker_{n}\Psi_{H}^{_{\tau}G.N})\leq \mathrm{J}_{_{\tau}G.N}(ker_{n}\Psi_{H}^{_{\tau}G.N}).$$
		Only verify $\phi(z_{i}e_{p-1})\leq \mathrm{J}_{_{\tau}G.N}(ker_{n}\Psi_{H}^{_{\tau}G.N})$ for all $1\leq i\leq n$. Assume $y_{i}=(y_{i,1},\cdots,y_{i,n})$ and $ z_{i}=(z_{i,1},\cdots,z_{i,n})$. Then $$\phi(z_{i}e_{p-1})=(\sum\limits_{j=1}^{n}z_{i,j}(a_{1}-1)^{p-1}y_{j,1},\cdots,\sum\limits_{j=1}^{n}z_{i,j}(a_{1}-1)^{p-1}y_{j,n}).$$
		And $\sum\limits_{j=1}^{n}z_{i,j}(a_{1}-1)^{p-1}y_{j,s}=\sum\limits_{j=1}^{n}z_{i,j}y_{j,s}(a_{1}-1)^{p-1}+\sum\limits_{j=1}^{n}z_{i,j}[(a_{1}-1)^{p-1},y_{j,s}]$ where $[(a_{1}-1)^{p-1},y_{j,s}]=(a_{1}-1)^{p-1}y_{j,s}-y_{j,s}(a_{1}-1)^{p-1}$. Since $[(a_{1}-1)^{p-1},y_{j,s}]\in ker\Psi_{H}^{_{\tau}G.N}$ and $z_{i}\in \mathrm{J}_{_{\tau}G.N}(\prod\limits^{n}\mathbb{F}_{p}(_{\tau}G.N))$, $\sum\limits_{j=1}^{n}z_{i,j}[(a_{1}-1)^{p-1},y_{j,s}]\in \mathrm{J}_{_{\tau}G.N}(ker\Psi_{H}^{_{\tau}G.N})$. Then $\sum\limits_{j=1}^{n}z_{i,j}y_{j,s}(a_{1}-1)^{p-1}\in ker\Psi_{H}^{_{\tau}G.N}.$ From Lemma \ref{xo}, $ \sum\limits_{j=1}^{n}z_{i,j}y_{j,s}(a_{1}-1)^{p-1}\in \mathrm{J}_{_{\tau}G.N}(ker\Psi_{H}^{_{\tau}G.N})$. Then $\phi(z_{i}e_{p-1})\in \mathrm{J}_{_{\tau}G.N}(ker_{n}\Psi_{H}^{_{\tau}G.N})$.

		The proof is completed.
	\end{pf}
	
	{\noindent\bf Proof of Theorem \ref{qq}}.
	
	\begin{pf}
		From Lemma \ref{ut}, $Q$ can be embedded in $\prod\limits^{n}\mathbb{F}_{p}(G) $ as $G$-module. From Theorem \ref{aa} and Theorem \ref{ggg}, the proof is completed.
	\end{pf}
	
	\begin{lem}\label{kj}
		Let $Q$ be a $1-G$-module. If  $\mathrm{H}^{1}(G,Q)\cong\mathrm{H}^{1}(_{\tau}G.\mathrm{C}_{p},Q)$ for $\tau\in Z^{2}(G,\mathrm{C}_{p})$, $$\mathrm{H}^{2}(_{\tau}G.\mathrm{C}_{p},Q)\cong \mathrm{H}^{2}(G,Q)\cong \mathbb{F}_{p}$$ and $Q$ is a cyclic $G$-module.
	\end{lem}
	\begin{pf}
		Set $Q_{1}$ is $G$-submodule of $\mathbb{F}_{p}(G)$ that $Q\cong Q_{1}$ as $G$-module. Denote $Q_{2}=\Psi_{_{\tau}G.\mathrm{C}_{p}}^{G}(Q_{1})$. Without loss of generality, suppose $Q_{1},Q_{2}$ are two right $G$-modules.
		
		Since $\mathrm{H}^{1}(_{\tau}G.\mathrm{C}_{p},Q)\cong\mathrm{F}_{p}$, set $x\in\mathbb{F}_{p}(_{\tau}G.\mathrm{C}_{p})$ that $\mathcal{L}_{_{\tau}G.\mathrm{C}_{p}}(Q_{2})=\mathbb{F}_{p}(_{\tau}G.\mathrm{C}_{p})x$.  From Lemma \ref{ii}, $ Ann_{L}(x)=\Psi_{_{\tau}G.\mathrm{C}_{p}}^{G}(Ann_{L}(\Psi_{G}^{_{\tau}G.\mathrm{C}_{p}}(x)))$.
		
		Set $ a\in\, _{\tau}G.\mathrm{C}_{p}\backslash \{1\}$ that $a-1\in  ker\Psi_{G}^{_{\tau}G.\mathrm{C}_{p}}$. From Lemma \ref{dd}, $ker\Psi_{G}^{_{\tau}G.\mathrm{C}_{p}} $ is left (right) cyclic $G$-module generated by $a-1$. Since $ker\Psi_{G}^{_{\tau}G.\mathrm{C}_{p}}\leq \mathcal{L}_{_{\tau}G.\mathrm{C}_{p}}(Q_{2})$, there exists $y\in \mathbb{F}_{p}(G)\backslash ker\Psi_{G}^{_{\tau}G.\mathrm{C}_{p}}$ such that $ yx=a-1$.
		
		Next, prove $Q_{2}$ is right $G$-module generated by $\Psi_{_{\tau}G.\mathrm{C}_{p}}^{_{\tau}G.\mathrm{C}_{p}}(y)$.
		Notice $a\in Z(_{\tau}G.\mathrm{C}_{p})$. Thus $xyx-yxx=0$. And $(xy-yx)x=0$. This imply $$ xy-yx\in Ann_{L}(x).$$ Notice $Ann_{L}(\Psi_{G}^{_{\tau}G.\mathrm{C}_{p}}(x))\leq\mathrm{J}_{G}(\mathbb{F}_{p}(G)) $. Then $xy-yx\in\mathrm{J}_{_{\tau}G.\mathrm{C}_{p}}(\mathrm{Im}\Psi_{_{\tau}G.\mathrm{C}_{p}}^{G} ) $. Thus $xy=y_{1}(a-1)$ where $y_{1}\in \mathbb{F}_{p}(_{\tau}G.\mathrm{C}_{p})\backslash \mathrm{J}_{_{\tau}G.\mathrm{C}_{p}}(\mathbb{F}_{p}(_{\tau}G.\mathrm{C}_{p}))$. This imply $ker\Psi_{G}^{_{\tau}G.\mathrm{C}_{p}}\leq x\mathrm{F}_{p}(_{\tau}G.\mathrm{C}_{p})$. Then $\mathcal{R}_{_{\tau}G.\mathrm{C}_{p}}(Ann_{L}(x))=x\mathrm{F}_{p}(G)$. From Lemma \ref{ii} and Remark \ref{rt}, $Q_{2}=Ann_{R}(x)$ is right cyclic $_{\tau}G.\mathrm{C}_{p}$-module. 
		
		Notice $x\mathbb{F}_{p}(_{\tau}G.\mathrm{C}_{p})\cong \mathbb{F}_{p}(_{\tau}G.\mathrm{C}_{p})/Q_{2}$ as right $_{\tau}G.\mathrm{C}_{p} $- module. Since $$\mathrm{H}^{2}(_{\tau}G.\mathrm{C}_{p},\mathbb{F}_{p}(_{\tau}G.\mathrm{C}_{p}))=0,$$ $\mathrm{H}^{2}(_{\tau}G.\mathrm{C}_{p},Q_{2})\cong\mathrm{H}^{1}(_{\tau}G.\mathrm{C}_{p},x\mathbb{F}_{p}(_{\tau}G.\mathrm{C}_{p}))$. Only verify $\mathcal{L}_{_{\tau}G.\mathrm{C}_{p}}(x\mathbb{F}_{p}(_{\tau}G.\mathrm{C}_{p}))$ is a cyclic left $_{\tau}G.\mathrm{C}_{p} $-module. Notice $\mathcal{L}_{_{\tau}G.\mathrm{C}_{p}}(x\mathbb{F}_{p}(_{\tau}G.\mathrm{C}_{p}))=Ann_{L}(x)$. From Lemma \ref{ii}, $\mathcal{L}_{_{\tau}G.\mathrm{C}_{p}}(x\mathbb{F}_{p}(_{\tau}G.\mathrm{C}_{p})) $ is cyclic. It's similar to prove $\mathcal{L}_{G}(x\mathbb{F}_{p}(_{\tau}G.\mathrm{C}_{p})/ker\Psi_{G}^{_{\tau}G.C_{p}})$ is cyclic module.
		
		The proof is completed.	 
	\end{pf}
	
	\begin{rem}
		Let $Q$ is module that described in the above Lemma. It's similar to prove $\mathrm{H}^{n}(_{\tau}G.\mathrm{C}_{p},Q)\cong H^{n}(G,Q)\cong \mathbb{F}_{p}$.
	\end{rem}
	
	\begin{thm}\label{dp}
		Let $Q_{1}$ be a $1-_{\tau}G.N$-module where $N$ is a $\mathbb{F}_{p}(G)$-module that $|N|=p^{2}$. Then $\mathrm{dim}_{\mathbb{F}_{p}}Q\geq p\mathrm{dim}_{\mathbb{F}_{p}}Q^{N}$.
	\end{thm}
	\begin{pf} Without loss of generality, suppose $Q$ is right $_{\tau}G.N$-module.
		Denote $Q_{1}$ as right $_{\tau}G.N$-module of $\mathbb{F}_{p}(_{\tau}G.N)$ that $Q\cong Q_{1}$ as right $_{\tau}G.N$-module. Assume $N=\langle a_{1}\rangle\times\langle a_{2}\rangle $. 
		
		When $\mathrm{H}^{1}(_{\tau}G.N,Q)=0$, $Q\cong \mathbb{F}_{p}(_{\tau}G.N)$ as right $_{\tau}G.N$ module. Then $Q^{N}\cong \mathbb{F}_{p}(G)$ as $G$- module. Thus $\mathrm{dim}_{\mathbb{F}_{p}}Q= p^{2}\mathrm{dim}_{\mathbb{F}_{p}}Q^{N}$. 
		
		When $\mathrm{H}^{1}(_{\tau}G.N,Q_{1})\cong\mathbb{F}_{p}$, $\mathcal{L}_{_{\tau}G.N}(Q)$ is a left cyclic $_{\tau}G.N$-module. Assume $ \mathcal{L}_{_{\tau}G.N}(Q_{1})=\mathbb{F}_{p}(_{\tau}G.N)y$.  Denote $I_{m}$ as left $_{\tau}G.N$-module similar to definition of $I_{n,m}$ in Lemma \ref{dd}. From Remark \ref{xp}, $I_{p-1}/I_{p}$ is isomorphic to direct sum of $p$ copies of $\mathbb{F}_{p}(G)$ as left(right) $G$-module. Set \begin{align*}
			\psi: \mathbb{F}_{p}(_{\tau}G.N)&\rightarrow \mathbb{F}_{p}(_{\tau}G.N)/I_{p}\\
			x&\mapsto xy.
		\end{align*} If $ y\in ker\Psi_{G}^{_{\tau}G.N}$, $I_{p-1}\leq ker \psi$. Then $\frac{|\mathbb{F}_{p}(_{\tau}G.N)/I_{p}|}{|\mathrm{Im}\psi|}\geq |\mathbb{F}_{p}(G)|^{p}$. And $\mathrm{dim}_{\mathbb{F}_{p}}(\mathcal{L}_{_{\tau}G.N}(Q_{1})+I_{p})\leq \mathrm{dim}_{\mathbb{F}_{p}}(\mathbb{F}_{p}(_{\tau}G.N)/I_{p})-p|G|$. Thus $\mathrm{dim}_{\mathbb{F}_{p}}Q_{1}\geq p|G|=p\mathrm{dim}_{\mathbb{F}_{p}}Q_{1}^{N}$. 
		
		Next, suppose $ y\notin ker\Psi_{G}^{_{\tau}G.N}$. Denote $A$ as $\mathbb{F}_{p}$ linear space over basis $CR(G,\,_{\tau}G.N)$ defined in Lemma \ref{xo}. And $B_{i}=\{(a_{1}-1)^{i-1}(a_{2}-1)^{p-i}x|x\in A\}$. From Lemma \ref{xo}, $I_{p-1}=B_{1}\oplus\cdots\oplus B_{p}\oplus I_{p}$. Set $ Y_{1}=\{x\in A|xy\in I_{1}\}$. And $ A_{1}$ is left $_{\tau}G.N$-module generated by $y,ker\Psi_{G}^{_{\tau}G.N}$. Then $A_{1}=\mathcal{L}_{_{\tau}G.N}(Q_{1}^{N})$. Notice  $\Psi_{G}^{_{\tau}G.N}(Y_{1})=\Psi_{G}^{_{\tau}G.N}(A_{1})=\{x\in \mathbb{F}_{p}(G)|x\Psi_{G}^{_{\tau}G.N}(y)=0\}$.  Thus $|Y_{1}|=|Q_{1}^{N}|$. Set $X_{i}=\{(a_{1}-1)^{i-1}(a_{2}-1)^{p-i}x|x\in Y_{1}\}$. For $\sum\limits_{i=1}^{p}(a_{1}-1)^{i-1}(a_{2}-1)^{p-i}x_{i}\in B_{1}\oplus\cdots\oplus B_{p}$ where $x_{i}\in A$ for all $1\leq i\leq n$, $ \sum\limits_{i=1}^{p}(a_{1}-1)^{i-1}(a_{2}-1)^{p-i}x_{i}\in ker\psi$ if and only if $x_{i}\in Y_{1}$ for all $1\leq i\leq n$. 
		Then $X_{1}\oplus\cdots\oplus X_{p}\leq ker\psi$. Thus $\mathrm{dim}_{\mathbb{F}_{p}}(\mathrm{Im}\psi)\leq \mathrm{dim}_{\mathbb{F}_{p}}(\mathbb{F}_{p}(_{\tau}G.N)/I_{p})-p\mathrm{dim}_{\mathbb{F}_{p}}(Y_{1})$. Since $|Y_{1}|=|Q_{1}^{N}|$, $ \mathrm{dim}_{\mathbb{F}_{p}}(\mathcal{L}_{_{\tau}G.N}(Q_{1}))\leq \mathrm{dim}_{\mathbb{F}_{p}}(\mathbb{F}_{p}(_{\tau}G.N))-p\mathrm{dim}_{\mathbb{F}_{p}}(Q_{1}^{N})$. Thus $\mathrm{dim}_{\mathbb{F}_{p}}(Q_{1})\geq p\mathrm{dim}_{\mathbb{F}_{p}}(Q_{1}^{N})$.
		
		The proof is completed.
	\end{pf}
	
	\begin{cor}\label{rty}
		Let $N$ be a $\mathbb{F}_{p}(G)$-module that $|N|=p^{2}$. $N_{1}$ is a $G$-submodule of $N$ that $|N_{1}|=p$. $Q$ is a $1-_{\tau}G.N/N_{1}$-module that $$\mathrm{H}^{1}(_{\tau}G.N/N_{1},Q)\cong\mathrm{H}^{1}(_{\tau}G.N,Q).$$
		Then $\mathrm{dim}_{\mathbb{F}_{p}}Q=p\mathrm{dim}_{\mathbb{F}_{p}}Q^{N}$. Specially, $Q^{N}$ is a cyclic $G$-module.
	\end{cor}
	\begin{pf}
		Without loss of generality, suppose $Q$ is right $_{\tau}G.N$-module. 
		
		From Theorem \ref{dp}, $\mathrm{dim}_{\mathbb{F}_{p}}Q\geq p\mathrm{dim}_{\mathbb{F}_{p}}Q^{N}$. Assume $N/N_{1}=\langle a \rangle$. Set $Q_{i}=\{x\in Q|x^{(a-1)^{i}}\in Q^{N}\}$. Then $Q_{p-1}=Q$. From $x^{a-1}=0$ if and only if $x\in Q^{N}$, $ |Q_{i+1}/Q_{i}|\leq |Q_{0}|=|Q^{N}|$ for all $i\in \mathbb{N}_{+}$. Since $Q_{i+1}=Q_{i}$ for all $i\geq p-1$ and $\mathrm{dim}_{\mathbb{F}_{p}}Q\geq p\mathrm{dim}_{\mathbb{F}_{p}}Q^{N}$, $|Q_{i+1}/Q_{i}|=|Q^{N}|$ for all $ 0\leq i\leq p-2$. Then $\mathrm{dim}_{\mathbb{F}_{p}}Q=p\mathrm{dim}_{\mathbb{F}_{p}}Q^{N}$.
		
		Set
		\begin{align*}
			\psi: Q&\rightarrow Q^{N}\\
			x&\mapsto x^{(a-1)^{p-1}}.
		\end{align*}
		Since $a\in Z(_{\tau}G.N/N_{1})$, $\psi$ is a $_{\tau}G.N/N_{1}$-module epimorphism from $Q$ to $Q^{N}$.  From Lemma \ref{kj}, $Q^{N}$ is cyclic $_{\tau}G.N$-module. Since $N$ acts on $Q^{N}$ trivially, $Q^{N}$ is a cyclic $G$-module.
		
		The proof is completed.
	\end{pf}

	\begin{lem}\label{xu}
		Let $Q$ be a $\mathbb{F}_{p}(G)$-module. $N$ is a minimal normal subgroup of $G$ that $ |Q^{N}|^{p}=|Q|$. Then $Q$ is direct sum of $n$ copies of $ \mathbb{F}_{p}(N)$ as $N$-module where $|Q^{N}|=p^{n}$.  And $\mathrm{H}^{1}(N,Q)=0$.
	\end{lem}
	\begin{pf}Set $a\in N\backslash \{1\}$. Denote $Q_{i}=\{x^{(a-1)^{i}}x\in Q\}$.
		Since $\mathrm{dim}_{\mathbb{F}_{p}}(Q)=\mathrm{dim}_{\mathbb{F}_{p}}(Q^{N})$, $|Q_{i}/Q_{i+1}|=|Q_{p-1}|$ for all $0\leq i\leq p-1$.  From $Q_{p-1}=Q^{N}$, $Q$ is direct sum of $n$ copies of $ \mathbb{F}_{p}(N)$ as $N$-module where $|Q^{N}|=p^{n}$.
		
		The proof is completed.
	\end{pf}
	
	The following Theorems will discuss the finite $p$-group with non-cyclic elementary abelian normal $p$-group. And the finite $p$-group with unique elementary abelian normal $p$-group will be dicussed in next section. 
	
	\begin{thm}\label{io}
		Let $G$ be a non-cyclic finite $p$-group. $N$ is a minimal normal subgroup of $G$. $Q$ is a $1-G/N$-module that $\mathrm{H}^{1}(G/N,Q)=\mathrm{H}^{1}(G,Q)$. Denote $\pi_{\tau}$ as surjective group homomorphism from $_{\tau}G.Q/\mathrm{J}_{G}(Q)$ to $G$ where $\tau\in Z^{2}(G,Q)$.  Then 
		 when $G$ has a normal elementary abelian subgroup of order equal to $p^{2}$,  there exists a normal subgroup $T_{1}$ of $_{\tau}G.Q/\mathrm{J}_{G}(Q)$ such that $T_{1}\times ker\pi_{\tau}=\pi_{\tau}^{-1}(T)$ for a normal subgroup $T$ of $G$ satisfy $N\leq T,T\cong \mathrm{C}_{p}^{2}$. Sepcially, if $T\leq Z(G)$, $ \pi_{\tau}^{-1}(T)\leq \Omega_{1}(Z(_{\tau}G.Q/J_{G}(Q)))$.
	\end{thm}
	\begin{pf}Since $Q$ is a $1-G/N$ module that $\mathrm{H}^{1}(G/N,Q)=\mathrm{H}^{1}(G,Q)$, $Q$ is a $1-G$-module. From Lemma \ref{kj}, $\mathrm{H}^{2}(G,Q)\cong\mathbb{F}_{p}$. When $\tau\in B^{2}(G,Q)$, $_{\tau}G.Q\cong Q\rtimes G$. The result is obvious. Next, suppose $\tau\in Z^{2}(G,Q)\backslash B^{2}(G,Q)$. Since $\mathrm{H}^{2}(G,Q)\cong\mathbb{F}_{p}$, $_{\tau}G.Q\cong\, _{\d}G.Q$ for all $\d\in Z^{2}(G,Q)\backslash B^{2}(G,Q)$. Thus only verify that the result is true for some $\d\in Z^{2}(G,Q)\backslash B^{2}(G,Q)$.
		
		Next, suppose $G$ has a normal elementary abelian subgroup which is isomorphic to $\mathrm{C}_{p}^{2}$. 
		
		There exists a normal subgroup $M$ of $G$ such that $M\cong\mathrm{C}_{p}^{2}$ and $N\leq M$. From Corllary \ref{rty}, $\mathrm{dim}_{\mathbb{F}_{p}}(Q)=p\mathrm{dim}_{\mathbb{F}_{p}}(Q^{M}) $. Set $\tau\in Z^{2}(G,Q)\backslash B^{2}(G,Q)$.  And
		$\pi$ is natural surjective group homomorphism from $_{\tau}G.Q$ to $ G$. Denote $N_{1}=\pi^{-1}(N),M_{1}=\pi^{-1}(M)$. Since $ N$ acts on $Q$ trivially, $N_{1}$ is abelian. Set $x\in M_{1}\backslash N_{1}$. When $M\nleqslant Z(G)$, there exists $g\in _{\tau}G.Q$ such that $[x,g]\in N_{1}\backslash Q$. And \begin{align*}[x^{p},g,y]&=[\prod\limits_{i=1}^{p}[x,g,\overbrace{x,\cdots,x}^{i-1}]^{p\choose i},y]\\&=\prod\limits_{i=1}^{p}[x,g,\overbrace{x,\cdots,x}^{i-1},y]^{p\choose i}\\&=[x,g,\overbrace{x,\cdots,x}^{p-1},y]
		\end{align*} for all $y\in _{\tau}G.Q$.
		Since $[x,g,y]\in N_{1}$, $[x^{p},g,y]=[x,g,y,\overbrace{x,\cdots,x}^{p-1}]=[x,g,y]^{(x-1)^{p-1}}$. Notice $x^{p}\in Q^{M}$. Then $[x^{p},g,y]\in \mathrm{J}_{G}(Q^{M}) $. Thus $[x,g,y]\in \mathrm{J}_{G}(Q)$. It's similar to prove $[x,z]\in \mathrm{J}_{G}(Q)$ for all $z\in _{\tau}G.Q$ when $M\leq Z(G)$. Then there exists a normal subgroup $T_{1}$ of $_{\tau}G.Q/\mathrm{J}_{G}(Q)$ such that $\pi_{\tau}^{-1}(M)=ker\pi_{_{\tau}}\times T_{1}$.
		
		The proof is completed.
	\end{pf}

	\begin{lem}\label{jx}
		Let $G$ be an non-cyclic abelian $p$-group. And $N$ is a cyclic subgroup of $G$. $Q$ is a $1-G$-module that $\mathrm{H}^{1}(G/N,  Q^{N})\cong\mathrm{H}^{1}(G,Q^{N})$. Then $ d(_{\tau}G.Q)=d(G)+1$ and $d(_{\d}(G/N).Q)=d(G)+1$ for all $\tau\in Z^{2}(G,Q),\d\in Z^{2}(G/N,Q)$. 
	\end{lem}
	\begin{pf}
		When $\tau\in B^{2}(G,Q)$, $_{\tau}G.Q\cong Q\rtimes G$. Then $d(_{\tau}G.Q)=d(G)+1$.

		From Lemma \ref{kj}, $\mathrm{H}^{2}(G,Q)\cong \mathrm{H}^{2}(G/N,Q)\cong \mathbb{F}_{p}$. Then the groups induced by all deriveions in $ Z^{2}(G,M)\backslash B^{2}(G,M)$  are isomorphic. We only prove there exists $\tau\in Z^{2}(G,Q)\backslash B^{2}(G,Q)$ such that $d(_{\tau}G.Q)=d(G)+1$. The other cases prove similarly.
		 
		Since  $G$ is abelian, there exists a abelian group $H$ and a surjective group homomorphism $\psi$ from $H$ to $G$ such that $ ker\psi\cong\mathrm{C}_{p}$ and $ \psi^{-1}(N)$ is cyclic. Then there exists $\tau\in Z^{2}(G,Q)\backslash B^{2}(G,Q)$ such that $ \tau(x,y)\in Q^{G}$ for all $x,y\in G$. This imply $d(_{\tau}G.Q)=d(G)+1$.
		
		The proof is completed.
	\end{pf}
	
	Denote $N$ as a normal subgroup of $G$. Set $A$ is a $G/N$-module. And \begin{align*} \pi:Z^{1}(G/N,A^{N})&\rightarrow Z^{1}(G,A)\\
		\tau&\mapsto \d
	\end{align*} where $\d(gN_{1})=\tau(gN)$ for all $g\in G$.
	$\pi$ is injective. In the following discussion, say $ \tau\in Z^{1}(G/N,A)\backslash Z^{1}(G,A^{N})$ means that $ \tau\in Z^{1}(G/N_{1},A)\backslash\mathrm{Im}\pi$.
	
	\begin{lem}\label{px}
		Let $N_{1},N_{2}$ be two distinct minimal normal subgroups of $G$. And $Q$ is a $1-G/N_{1}$ module that $\mathrm{dim}_{\mathbb{F}_{p}}(Q)=p\mathrm{dim}_{\mathbb{F}_{p}}Q^{N_{1}\times N_{2}}$. Then $\mathrm{H}^{1}(G,Q^{N_{1}\times N_{2}})\cong\mathrm{H}^{1}(G/N_{1},Q^{N_{1}\times N_{2}})$ if and only if $\mathrm{H}^{1}(G,Q)\cong \mathrm{H}^{1}(G/N_{1},Q)$. 
	\end{lem}
	\begin{pf}
		Set $Q_{i}=\{x^{(a-1)^{i}}|x\in Q\}$ where $N_{2}=\langle a\rangle$.
		Similar to the proof of Corllary \ref{rty}, $|Q_{i}/Q_{i+1}|=|Q^{N_{1}\times N_{2}}|$.
		
		When $\mathrm{H}^{1}(G,Q)\cong \mathrm{H}^{1}(G/N_{1},Q)$, if $\mathrm{H}^{1}(G,Q^{N_{1}\times N_{2}})\geq\mathrm{H}^{1}(G/N_{1},Q^{N_{1}\times N_{2}})\oplus\mathbb{F}_{p}$, there exists $\tau\in Z^{1}(G,Q^{N_{1}\times N_{2}})\backslash Z^{1}(G/N_{1},Q^{N_{1}\times N_{2}})$. Then $\tau\in Z^{1}(G,Q)\backslash Z^{1}(G/N_{1},Q)$. It's a contradiction.
		
		Next, suppose $\mathrm{H}^{1}(G,Q^{N_{1}\times N_{2}})\cong\mathrm{H}^{1}(G/N_{1},Q^{N_{1}\times N_{2}})$. If $\mathrm{H}^{1}(G,Q)\geq \mathrm{H}^{1}(G/N_{1},Q)\oplus\mathbb{F}_{p}$, there exists $\tau\in Z^{1}(G,Q)\backslash Z^{1}(G/N_{1},Q)$. Since $a^{p}=1$, $\tau(a)^{(a-1)^{p-1}}=0$. Then $\tau(a)\in Q_{1}$. Thus there exists $\d\in B^{1}(G,Q)$ such that $(\tau+\d)(a)=0$. From $[a,x]=1$ for all $x\in G$, $0=(\tau+\d)([a,x])=(\tau+\d)(x)^{1-a}$. Then $ (\tau+\d)(x)\in Q^{N_{1}\times N_{2}}$. This imply $(\tau+\d)\in Z^{1}(G,Q^{N_{1}\times N_{2}})$. Then $\mathrm{H}^{1}(G,Q^{N_{1}\times N_{2}})\geq\mathrm{H}^{1}(G/N_{1},Q^{N_{1}\times N_{2}})\oplus \mathbb{F}_{p}$. It's a contradiction.
		
		The proof is completed.
	\end{pf}
	
	\begin{lem}\label{pl}
		Let $N$ be a minimal normal subgroup of $G$. And $N_{1}$ is a normal subgroup of $G$ that satisfying $ N\leq N_{1}\cong C_{p}^{2}, N_{1}\cap Z(G)=N$. $Q$ is a $1-G/N$-module that $dim_{\mathbb{F}_{p}}Q=pdim_{\mathbb{F}_{p}}Q^{N_{1}}\geq 3$. Then $\mathrm{H}^{1}(G,Q)\cong \mathrm{H}^{1}(G/N,Q)$ if and only if
		$$\mathrm{H}^{1}(G,J_{G}(Q))\cong \mathrm{H}^{1}(G,J_{G}(Q)).$$
		\end{lem}
	\begin{pf}
		When $\mathrm{H}^{1}(G,Q)\cong \mathrm{H}^{1}(G/N,Q)$, it's easy to prove $\mathrm{H}^{1}(G,J_{G}(Q))\cong \mathrm{H}^{1}(G,J_{G}(Q)).$ Next, we prove the sufficiency.
		
	From Corllary\ref{8.0},	assume $\mathrm{H}^{1}(G,Q)\cong \mathrm{H}^{1}(G/N,Q)\oplus\mathbb{F}_{p}$. Set $$\tau\in Z^{1}(G,Q)\backslash Z^{1}(G/N,Q).$$ Suppose $Q_{i}=\{x^{(a-1)^{i}}|x\in Q\}$ where $0\leq i\leq p-1,a\in N_{1}\backslash N$. Since $dim_{\mathbb{F}_{p}}Q=pdim_{\mathbb{F}_{p}}Q^{N_{1}}$, it's easy to prove $ Q_{i}/Q_{i+1}\cong Q_{p-1}=Q^{N_{1}}$. Notice $0=\tau(a_{1}^{p})=\tau(a_{1})^{(a-1)^{p-1}}$. Then $\tau(a_{1})\in Q_{1}$. There exists $\d\in B^{1}(G,Q)$ such that $(\tau+\d)(a_{1})=0$. Since the sylow $p$-subgroup of  $Aut(C_{p}^{2})$ is isomorphic to $C_{p}$, there exists maximal subgroup $T$ of $G$ such that $[N_{1},T]=0$. Then $0=(\tau+\d)([a_{1},g])=(\tau+\d)(g)^{1-a_{1}}$ for all $g\in T$. This imply $(\tau+\d)(g)\in Q^{N_{1}}$ for all $g\in T$. Set $h\in G\backslash T$. Then $$(\tau+\d)([h,a_{1}])=(\tau+\d)(h)^{a_{1}-1}\in Q^{G}.$$
		From $(\tau+\d)(g)\in Q^{N_{1}}$ for all $g\in T$, then $ (\tau+\d)(h)^{x-1}\in Q^{N_{1}}$ for all $x\in G$. When $p\geq 3$, from $dim_{\mathbb{F}_{p}}Q=pdim_{\mathbb{F}_{p}}Q^{N_{1}}$, $(\tau+\d)(x)\in Q_{1}\in J_{G}(Q)$. When $p=2$, since $dim_{\mathbb{F}_{p}}Q\geq 3$, $ dim_{\mathbb{F}_{p}}Q^{N_{1}}\geq 2$. From $Q$ is cyclic module, $(\tau+\d)(x)\in  J_{G}(Q)$.
	Then $(\tau+\d)\in Z^{1}(G,J_{G}(M))\backslash Z^{1}(G/N,J_{G}(M))$.
	
		The proof is completed.	
	\end{pf}

	\begin{thm}\label{du}
		Let $G$ be a finite $p$-group with non-cyclic elementary abelian normal $p$-group.
		Let $N$ be a cyclic normal subgroup of $G$. $Q$ is a $1-G/N$-module that $\mathrm{H}^{1}(G/N,Q)\cong\mathrm{H}^{1}(G,Q)$. Then 
		\begin{enumerate}
			\item[\rm(1)] if $Z(G)$ is non-cyclic, $$
			\mathrm{H}^{1}(_{\tau}G.Q/Q,J_{G}(M))\cong\mathrm{H}^{1}(G,J_{G}(Q))\oplus\mathbb{F}_{p} $$ for any $ \tau\in Z^{2}(G,Q)$.
			\item[\rm(2)] if $d(Z(G))=1$  and $ dim_{\mathbb{F}_{p}}Q\geq 3$, $$
			\mathrm{H}^{1}(_{\tau}G.Q/Q,J_{G}(Q))\cong\mathrm{H}^{1}(G,J_{G}(Q))\oplus\mathbb{F}_{p} $$ for any $ \tau\in Z^{2}(G,Q)$.
		\end{enumerate}
	\end{thm}
	\begin{pf}
		Assume $K={}_{\tau}G.Q/J_{G}(Q)$.
		From Theorem \ref{ggg}, $$\mathrm{H}^{1}(K,M)\cong\mathrm{H}^{1}(G,M)\oplus\mathbb{F}_{p}.$$ Then there exists minimal normal subgroup $N_{1}$ of $K$ such that $ \d(x)=0$ for all $\d\in Z^{1}(K,M),x\in N_{1}$. Set $\pi_{\tau}$ is natural surjective homomorphism from $K$ to $G$. Since $\mathrm{H}^{1}(G/N,Q)\cong\mathrm{H}^{1}(G,Q)$ and $\mathrm{H}^{1}(K,M)\cong\mathrm{H}^{1}(G,M)\oplus\mathbb{F}_{p}$, $ N_{1}\cap ker\pi_{\tau}=1,N_{1}\leq \pi_{\tau}^{-1}(N)$. Then $$\mathrm{H}^{1}(K/N_{1},M)\cong \mathrm{H}^{1}(G/N,M)\oplus\mathbb{F}_{p}.$$
		
		When $Z(G)$ is non-cyclic, there exists a normal subgroup $N_{2}$ of $G$ such that $ N\cap N_{2}\cong C_{p},N_{2}\cong C_{p}^{2},N_{2}\leq Z(G)$. From Corllary \ref{rty}, $dim_{\mathbb{F}_{p}}Q=pdim_{\mathbb{F}_{p}}Q^{N_{2}}$. From Lemma \ref{px}, $$ \mathrm{H}^{1}(K/N_{1},J_{G}(M))\cong \mathrm{H}^{1}(K/\pi_{\tau}^{-1}(N_{2}),J_{G}(M))\oplus\mathbb{F}_{p}.$$
		Then $\mathrm{H}^{1}(_{\tau}G.Q/Q,J_{G}(M))\cong\mathrm{H}^{1}(G,J_{G}(Q))\oplus\mathbb{F}_{p}$.
		
		When $d(Z(G))=1$  and $ dim_{\mathbb{F}_{p}}Q\geq 3$, there exists a normal subgroup $N_{3}$ of $G$ such that $ N\cap N_{3}\cong C_{p},N_{3}\cong C_{p}^{2},N_{3}\leq Z(G)$. From 
		Corllary \ref{rty} and Lemma \ref{pl}, $$\mathrm{H}^{1}(_{\tau}G.Q/Q,J_{G}(M))\cong\mathrm{H}^{1}(G,J_{G}(Q))\oplus\mathbb{F}_{p}.$$
		
		The proof is completed.
	\end{pf}
	
	In the last of the section, we give a Lemma that will be used in next section. The lemma discusses a special case of finite $p$-group with non-cyclic elementary abelian normal $p$-group. 
	
	\begin{lem}\label{kko}
			Let $G$ be a finite $p$-group with non-cyclic elementary abelian normal $p$-group. $N$ is a cyclic normal subgroup of $G$. $Q$ is a $1-G$-module that satisfying $\mathrm{H}^{1}(G,Q)\cong \mathrm{H}^{1}(G/N,Q),|Q|=p^{2}$. Then $G$ is isomorphic to one of the following groups 
			$$ M_{2^{n}},C_{2}\times C_{2^{m-1}}$$
			where $n\geq 4,m\geq 2$. And $ Q$ is isomorphic to $\mathbb{F}_{p}(G/T)$ where $T$ is a cyclic maximal subgroup of $G$  
	\end{lem}
\begin{pf}
	Firstly, we assert $d(G)=2$. Otherwise, $G/\Phi(G)\cong C_{p}^{s}$ where $s\geq 3$. From $M$ is a $1-G/N$-module and $$ \mathrm{H}^{1}(G,Q^{G})\cong Hom_{\mathbb{Z}}(G,Q^{G})\cong \mathbb{F}_{p}^{s},$$ $|M|\geq p^{3}$. It's a contradiction.  Set $K$ is a normal subgroup of $G$ that satisfying $N\cap K\cong C_{p}, K\cong C_{p}^{2}$. We assert $K\nleqslant \Phi(G)$. Otherwise, $dim_{\mathbb{F}_{p}}Q^{\Phi(G)}\geq 2$ from $ d(G)=2$. Since $K\leq \Phi(G)$, $dim_{\mathbb{F}_{p}}Q\geq pdim_{\mathbb{F}_{p}}Q^{\Phi(G)}\geq 2p>2$. It's  conratdictory to $|M|=p^{2}$. From Corllary \ref{rty}, $Q^{K}\cong \mathbb{F}_{p}$ and $p=2$. Then $G/K$ is cyclic. Thus $G$ is isomorphic to one of the following groups:
	$$M_{2^{n}},C_{2}\times C_{2^{m-1}}$$  
	where $n\geq 3,m\geq 2$. 
	
	Since sylow $2$-subgroup of $Aut(C_{2}^{2})$ is isomorphic to $C_{2}$, there exists a maximal subgroup $T$ of $G$ such that $ T$ acts on $M$ trivially. Next, we prove $T$ is cyclic. Otherwise, there exists a normal subgroup $T_{1}$ of $G$ such that $T_{1}\leq T,T_{1}\cong C_{2}^{2}$. From Theorem \ref{dp}, $$dim_{\mathbb{F}_{p}}Q\geq dim_{\mathbb{F}_{p}}Q^{T_{1}}\geq  dim_{\mathbb{F}_{p}}Q^{T}>2 .$$
	It's a contrdiction. Then $T$ is cyclic.
	
	The proof is completed.
\end{pf}

	\section{Automorphism of Order $p$}

	\begin{defi}
		Let $G$ be a non-abelian finite $p$-group. Say a normal subgroup $N$ of $G$ is \textbf{special} subgroup if 
		\begin{align*}
			&C_{G}(N)/Z(N)\;is\;a \;cyclic;\; and\\
			&\mathrm{I}(C_{G}(N))\leq N.
		\end{align*}
		
	\end{defi}

	\begin{lem}\label{ty}
		Denote $N$ as a special subgroup of $G$. Suppose $N<C_{G}(N)N$. If $\mathrm{H}^{1}(G/N,\Omega_{1}(Z(N)))\geq \mathrm{H}^{1}(G/C_{G}(N)N,\Omega_{1}(Z(N)))\oplus\mathbb{F}_{p}$ , then $G$ has a non-inner automorphism of order $p$. 
	\end{lem}
	\begin{pf}
		Notice $C_{G}(C_{G}(N)N)\leq C_{G}(N)$. Then $\Omega_{1}(Z(C_{G}(N)N))\leq \mathrm{I}(C_{G}(N))\leq N$. And $\Omega_{1}(Z(C_{G}(N)N))=\Omega_{1}(Z(N))$. Denote $$\tau\in Z^{1}(G/N,\Omega_{1}(Z(N)))\backslash Z^{1}(G/C_{G}(N)N,\Omega_{1}(Z(N))).$$  Set $\psi(g)=g\tau(gN)$ for all $g\in G$. The order of $\psi$ is $p$. If $\psi$ is inner, from $\psi(g)=g$ for all $g\in N$, then exists $h\in \mathrm{I}(C_{G}(N))$ such that $g^{h}=\psi(g)$ for all $g\in G$. Notcie $\mathrm{I}(C_{G}(N))\leq N$. Then $x^{h}=x$ for all $x\in C_{G}(N)N $. It's contradictory to $\tau\in Z^{1}(G/N,\Omega_{1}(Z(N)))\backslash Z^{1}(G/C_{G}(N)N,\Omega_{1}(Z(N)))$.
		
		The proof is completed.
	\end{pf}
	
	\begin{lem}\label{j}
		Let $A$ be a normal subgroup of $G$ that $\mathrm{I}(C_{G}(A))\leq A\leq\Phi(G)$. Assume $A_{1}$ is a normal subgroup of $G$ that $ |A_{1}/A|=p,Z(G)\leq Z(A_{1})$ and $$ \mathrm{H}^{1}(G/A_{1},\Omega_{1}(Z(A_{1})))\cong\mathrm{H}^{1}(G/A,\Omega_{1}(Z(A_{1}))).$$ Then $A_{1}\leq\Phi(G)$.
	\end{lem}
	\begin{pf}
		Notice $\mathrm{H}^{1}(G/\Phi(G),\Omega_{1}(Z(G)))\cong\mathrm{Hom}_{\mathbb{Z}}(G/\Phi(G),\Omega_{1}(Z(G)))$ and $$\mathrm{H}^{1}(G/A_{1}\Phi(G),\Omega_{1}(Z(G)))\cong\mathrm{Hom}_{\mathbb{Z}}(G/A_{1}\Phi(G),\Omega_{1}(Z(G))).$$
		If $A_{1}\nleqslant \Phi(G)$, $\mathrm{H}^{1}(G/\Phi(G),\Omega_{1}(Z(G)))\geq \mathrm{H}^{1}(G/A_{1}\Phi(G),\Omega_{1}(Z(G)))\oplus\mathbb{F}_{p}$. There exists $\tau\in Z^{1}(G/\Phi(G),\Omega_{1}(Z(G)))\backslash Z^{1}(G/A_{1}\Phi(G),\Omega_{1}(Z(G)))$. Since $A_{1}\nleqslant \Phi(G)$, $g\in A_{1}\Phi(G)\backslash\Phi(G)$ for all $g\in A_{1}\backslash A$. Set $\tau_{1}(gA)=\tau(g\Phi(G))$ for all $g\in G$. Then $\tau_{1}\in Z^{1}(G/,\Omega_{1}(Z(A_{1})))$. Notice $\tau_{1}(gA)\neq0$ for all $g\in A_{1}\backslash A$. Thus $\tau_{1}\notin Z^{1}(G/A_{1},\Omega_{1}(Z(A_{1})))$. It's  contradictory to $$\mathrm{H}^{1}(G/A_{1},\Omega_{1}(Z(A_{1})))\cong\mathrm{H}^{1}(G/A,\Omega_{1}(Z(A_{1}))).$$ The proof is completed.
	\end{pf}

	\begin{lem}\label{3.2}
		Assume $N$ is a normal subgroup of $G$ that $N\leq\Phi(G)$. If $$\mathrm{H}^{1}(G/N,\Omega_{1}(Z(N)))=1,$$ then $ \Omega_{1}(Z(N))$ is a proper subgroup of $N$.
	\end{lem}
	\begin{pf}
		If $ \Omega_{1}(Z(N))=N$, then $G$ is a extension of $G/N$ by $\Omega_{1}(Z(N))$. From main result in \cite{hoe}, $\mathrm{H}^{2}(G/N,\Omega_{1}(Z(N)))=1$. Then $ G\cong N\rtimes G/N$. Notice $d(N\rtimes G/N)>d(G/N)=d(G)$. It's contradictory to $ G\cong N\rtimes G/N$. Then $ \Omega_{1}(Z(N))$ is a proper subgroup of $N$.
	\end{pf}
	
	\begin{lem}\label{xi}
		Let $N$ be a special subgroup of $G$. Suppose all automorphisms induced by derivations in $Z^{1}(G/N,W)$ are inner where $W=\Omega_{1}(Z(N))$ and $n=d(Z(G))$. Then $C_{W}(A)$ is $n-G/A$-module for aribitrary normal subgroup $A$ of $G$ that $N\leq A$. 
	\end{lem}
	\begin{pf}
		From Theorem \ref{5.5}, $W$ is $n-G/N$-module. Notice $C_{G}(A)\leq C_{G}(N)$. Then $ A$ is a special subgroup of $G$. From Theorem 9.84 in \cite{rot}, $\mathrm{H}^{1}(G/A,C_{W}(A))\leq \mathbb{F}_{p}^{n}$.
		
		The proof is completed.
	\end{pf}

	\begin{lem}\label{yu}
		Let $N$ be a special subgroup of $G$. Set $A_{1},A_{2}$ are two normal subgroups of $G$ that $N\leq A_{2}<A_{1}$. Denote $W=\Omega_{1}(Z(N))$. Suppose all automorphisms induced by derivations in $Z^{1}(G/N,W)$ are inner.  If $\mathrm{H}^{1}(G/A_{2},C_{W}(A_{1}))\geq\mathrm{H}^{1}(G/N,W)\oplus\mathbb{F}_{p}$, then $$C_{W}(A_{1})<C_{W}(A_{2})$$.
	\end{lem}
	\begin{pf}
		From Theorem 9.84 in \cite{rot}, $ \mathrm{H}^{1}(G/A_{1},C_{W}(A_{1}))\leq \mathrm{H}^{1}(G/N,W)$. Set \begin{align*}
			\psi:Z^{1}(G/A_{2},C_{W}(A_{1}))&\rightarrow Z^{1}(G/N,W)\\
			\tau&\mapsto \tau_{1}
		\end{align*}
		where $\tau_{1}(x)=\tau(xA_{1})$. Since $\mathrm{H}^{1}(G/A_{2},C_{W}(A_{1}))\geq\mathrm{H}^{1}(G/N,W)\oplus\mathbb{F}_{p}$, $$ \psi^{-1}(B^{1}(G/N,W))>B^{1}(G/A_{2},C_{W}(A_{1})).$$ There exists $ h\in W\backslash C_{W}(A_{1})$ such that $\d_{h}\in \psi(Z^{1}(G/A_{2},A_{1}))$ where $\d_{h}(g)=g^{-1}g^{h}$. 
		
		The proof is completed.
	\end{pf}
	
	\begin{thm}\label{hh}
		Let $N$ be a self-centralizer normal subgroup of $G$ that $ N\leq\Phi(G)$. Denote $n=d(Z(G))$. If\; $\mathrm{H}^{1}(G/N,\Omega_{1}(Z(N)))\cong\mathbb{F}_{p}^{m}$ where $m< n$, then  either  $G$ has a non-inner automorphism of order $p$ or exists a normal subgroup $N_{1}$ of $G$ such that $C_{G}(N_{1})\leq N_{1}<N$.
	\end{thm}
	\begin{pf}
		Suppose all automorphisms induced by $Z^{1}(G/N,\Omega_{1}(Z(N)))$ are inner. From Lemma \ref{3.2}, $\Omega_{1}(Z(N))$ is a proper subgroup of $N$. Set $N_{1}$ is a normal subgroup of $G$ that $|N/N_{1}|=p$ and $ \Omega_{1}(Z(N))\leq N_{1}\leq N$. Since $\mathrm{H}^{1}(G/N,\Omega_{1}(Z(N)))\cong\mathbb{F}_{p}^{m}$ where $m< d(Z(G))$, from Theorem \ref{gg}, $\mathrm{H}^{1}(G/N_{1},\Omega_{1}(Z(N)))\geq \mathrm{H}^{1}(G/N,\Omega_{1}(Z(N)))\oplus\mathbb{F}_{p}^{n-m}$.
		
		Suppose $C_{G}(N_{1})\nleqslant N_{1}$. 
		Denote $W=\Omega_{1}(Z(N))$.
		
		If $N<C_{G}(N_{1})N$, set $A$ is a normal subgroup of $G$ that $ N\leq A\leq  C_{G}(N_{1})N$ and $|A/N|=p.$ From Lemma \ref{xi}, $C_{W}(A)$ is $n-G/A$-module. By Theorem \ref{gg}, $\mathrm{H}^{1}(G/N,C_{W}(A))\geq\mathbb{F}_{p}^{n}>\mathrm{H}^{1}(G/N,W)$. From Lemma \ref{yu}, $C_{W}(A)<W$. Since $ W\leq N_{1}$, $x^{h}= x$ for all $h\in W\backslash N,x\in (A\cap C_{G}(N_{1}))\backslash N$. It's a contradiction.
		
		If $ N=C_{G}(N_{1})N_{1}=N$, from Theorem \ref{gg}, $$\mathrm{H}^{1}(G/N_{1},\Omega_{1}(Z(N)))\geq \mathrm{H}^{1}(G/N,\Omega_{1}(Z(N)))\oplus\mathbb{F}_{p}^{n-m}.$$  There exists $\tau\in Z^{1}(G/N_{1},W)\backslash Z^{1}(G/N,W)$. Set $\phi$ is automorphism induced by $\tau$. Then $\phi(x)\neq x$ for all $x\in C_{G}(N_{1})\backslash N_{1}$. Since $|N/N_{1}|=p$, $C_{G}(N_{1})$ is abelian. Notice $\phi(x)=x$ for all $x\in N_{1}$. Thus $\phi$ is non-inner. 
		
		The proof is completed.
		
	\end{pf}

	\begin{lem}\label{ll}
		Let $N$ be a special subgroup of $G$. And $d(Z(G))=n$. Suppose $\Omega_{1}(Z(N))$ is a exactly $n-G/N$-module. And all automorphisms induced by derivations in $Z^{1}(G/N,\Omega_{1}(Z(N)))$ are inner. Set $N_{1}$ is a normal subgroup of $G$ that $\Omega_{1}(Z(N))Z(G)\leq N_{1}$. Then $ C_{G}(N_{1})N/N$ is cyclic.
	\end{lem}
	\begin{pf} 
		Suppose $C_{G}(N_{1})N/N $ is non-cyclic.
		Denote $W=\Omega_{1}(Z(NC_{G}(N)))$.

		Since $C_{G}(N_{1})N/N$ isn't cyclic, $\Phi(C_{G}(N_{1})N/N)$ is normal subgroup of $G/N$ that $ |(C_{G}(N_{1})N/N)/\Phi(C_{G}(N_{1})N/N)|\geq C_{p}^{2}$. Set $A$ is normal subgroup of $G$ that $A/N=\Phi(C_{G}(N_{1})N/N)$. Then $C_{G}(N_{1})N/A\geq C_{p}^{2}$. Thus exists a normal subgroup $B$ of $G$ such that $A\leq B\leq C_{G}(N_{1})N$ and $C_{G}(N_{1})N/B\cong C_{p}^{2}$. Denote $E=C_{G}(N_{1})N$. Claim $\mathrm{H}^{1}(G/B,C_{W}(E))>\mathbb{F}_{p}^{n}$. From Lemma \ref{xi}, $ C_{W}(E)$ is a $n-G/E$-module. From Theorem \ref{qq}, $\mathrm{H}^{1}(G/B,C_{W}(E))\geq \mathbb{F}_{p}^{n+1}$. By Lemma \ref{yu}, $C_{W}(E)<C_{W}(B)$. Since $W\leq N_{1}$, there exists $ h\in C_{W}(B)\backslash C_{W}(E)$ such that $[h,x]=1 $ for all $x\in (E\cap C_{G}(N_{1}))\backslash B$.
		It's contradictory to $h\in C_{W}(B)\backslash C_{W}(E)$.
		
		The proof is completed.
	\end{pf}
	
	\begin{lem}\label{qp}
		Let $N$ be a normal subgroup of $G$ that $\mathrm{I}(C_{G}(N))\leq N\leq \Phi(G)$.  If $NC_{G}(N)\nleqslant \Phi(G)$, there exists a non-inner automorphism of order $p$.
	\end{lem}
	\begin{pf}
		Since $NC_{G}(N)\nleqslant \Phi(G)$, there exists $h\in C_{G}(N)\backslash \Phi(G)$. Then there exists a maximal subgroup $M$ of $G$ such that $h\notin M$. Notice  $\mathrm{H}^{1}(G/M,\Omega_{1}(Z(G)))\cong\mathrm{Hom}_{\mathbb{Z}}(G/M,\Omega_{1}(Z(G)))$. Set $\tau\in Z^{1}(G/M,\Omega_{1}(Z(G)))\backslash B^{1}(G/M,\Omega_{1}(Z(G)))$. And $\psi(g)=g\tau(gM)$. Since $\Omega_{1}(Z(G))\leq N$, $\psi$ is order $p$. From $\mathrm{I}(C_{G}(N))\leq N$, $\psi$ is non-inner. Otherwise, set $x\in \mathrm{I}(C_{G}(N))$ that $ g^{x}=\psi(g)$ for all $g\in G$. From $\psi(h)\neq 1$, $ [x,h]\neq 1$. Notice $ h\in C_{G}(N)$. It's a contradiction.
		
		The proof is completed.
	\end{pf}

	\begin{lem}\label{kl}
		Let $N$ be a special subgroup of $G$.  Suppose $\Omega_{1}(Z(N))$ is a exactly $n-G/N$-module where $n=d(Z(G))$. And all automorphisms induced by derivations in $Z^{1}(G/N,\Omega_{1}(Z(N)))$ are inner. Then exists a normal subgroup $N_{1}$ of $G$ that $\Omega_{1}(Z(N))Z(G)\leq N_{1}< N,|N/N_{1}|=p$. Furthermore,  $ C_{G}(N)N/N_{1}$ is non-cyclic when $C_{G}(N)N/N$ is non-trivially . 
	\end{lem}
	\begin{pf}
		Since all automorphisms induced by derivations in $Z^{1}(G/N,\Omega_{1}(Z(N)))$ are inner, $ |\mathrm{I}(C_{G}(N))/\Omega_{1}(Z(N))Z(G)|\geq p$. Then exists a normal subgroup $N_{1}$ of $G$ that $\Omega_{1}(Z(N))Z(G)\leq N_{1}< N,|N/N_{1}|=p$.
		
		Next, suppose $C_{G}(N)/Z(N)$ is a non-trivial cyclic $p$-group. 
		
		Set $A$ is normal subgroup of $G$ that $A/N\cong\mathrm{C}_{p}$ and $ A\leq C_{G}(N)N$. Denote $a\in (A\cap C_{G}(N))\backslash N$. 
		
		Since $A/\Omega_{1}(Z(N))Z(G)$ is non-cyclic, there exists a normal subgroup $E$ of $G$ such that $\Omega_{1}(Z(N))Z(G)\leq E\leq N$ and $N/E\cong\mathrm{C}_{p}^{2}$. When $ A/E\cong C_{p}^{3}$, set $N_{1}$ is  a normal subgroup of $G$ that $E\leq N_{1}\leq N,|N/N_{1}|=p$. Then $A/N_{1}\cong C_{p}^{2}$.
		 When $A/E\cong C_{p^{2}}\times C_{p}$, since $\mho_{1}(C_{p^{2}}\times C_{p})\;char\; C_{p^{2}}\times C_{p}$ and $N/E=\Omega_{1}(A/E)$, there exists a normal subgroup $T$ of $G$ that $A/T\cong C_{p}^{2},T\leq N$.
		
		The proof is completed.
	\end{pf}
	
	\begin{lem}\label{qk}
		Let $N$ be a special subgroup of $G$. Suppose $\Omega_{1}(Z(N))$ is a exactly $n-G/N$-module where $n=d(Z(G))$. And all automorphisms induced by derivations in $Z^{1}(G/N,\Omega_{1}(Z(N)))$ are inner. Suppose $C_{G}(N)/Z(N)$ is a non-trivial cyclic $p$-group. Assume $N_{1}$ is a normal subgroup of $G$ that $\Omega_{1}(Z(N))Z(G)\leq N_{1}<N,|N/N_{1}|=p$ and $ C_{G}(N)N/N_{1}$ non-cyclic. If $ C_{G}(N_{1})N_{1}/N_{1}$ is non-cyclic, then $G$ has a non-inner automorphism of order $p$.
	\end{lem}
	\begin{pf}
		Denote $W=\Omega_{1}(Z(N)).$ Notice $ C_{G}(NC_{G}(N))\leq C_{G}(N)\leq NC_{G}(N)$. Then $\mathrm{I}(C_{G}(N))=\mathrm{I}(C_{G}(NC_{G}(N)))$ and $\Omega_{1}(Z(N))=\Omega_{1}(Z(NC_{G}(N)))$. From Lemma \ref{ty} and Theorem \ref{gg}, suppose $\mathrm{H}^{1}(G/NC_{G}(N),W)=\mathrm{H}^{1}(G/N,W)\cong\mathbb{F}_{p}^{n}$.

		From Lemma \ref{ll}, suppose $ C_{G}(N_{1})N/N$ is cyclic. Since $C_{G}(N_{1})N_{1}/N_{1}$ isn't cyclic, $N\leq N_{1}C_{G}(N_{1})$. Then $C_{G}(N_{1})N_{1}/N_{1}\cong C_{p^{m}}\times C_{p}$ where $ C_{G}(N_{1})N/N\cong C_{p^{m}}$. Set $D$ is normal subgroup that $N_{1}\leq D\leq C_{G}(N_{1})N_{1}$ and $ D/N_{1}\cong C_{p}^{2}$. From (2) in Theorem \ref{qq}, there exists $\a\in Z^{1}(G/N_{1},W)\backslash Z^{1}(G/D,W)$ such that $\a(g)\neq1$ for some $g\in N\backslash N_{1}$. Set $\rho(g)=g\a(gN_{1})$ for all $g\in G$. $\rho$ is a automorphism of order $p$ fixes $N_{1}$ not $N$. If $\rho$ is inner, then exists $ h\in C_{G}(N_{1})\backslash N$ such that $g^{h}=\rho(g)$ for all $g\in G$. 
		Since $\rho$ is order $p$, $h^{p}\in Z(G)\leq N_{1}$.
		And exists $a\in C_{G}(N)$ such that $ah\in C_{G}(N_{1})\cap N$. Notice $C_{G}(N_{1})\cap N $ is abelian. Then $[ah,g]=1$ for $g\in (C_{G}(N_{1})\cap N)\backslash N_{1}$. It's contradictory to $\rho(g)\neq1$ for some $g\in N\backslash N_{1}$. Then $\rho$ is non-inner.

		The proof is completed.
	\end{pf}

	\begin{thm}\label{xx}
		Let $N$ be a special subgroup of $G$ that $N/Z(G)\Omega_{1}(Z(N))$ is non-cyclic. Suppose $\Omega_{1}(Z(N))$ is a exactly $n-G/N$-module where $n=d(Z(G))$. And all automorphisms induced by derivations in $Z^{1}(G/N,\Omega_{1}(Z(N)))$ are inner. Suppose $C_{G}(N)/Z(N)$ is a non-trivial cyclic $p$-group. 
		Then one of following holds:
		\begin{enumerate}
			\item[\rm(1)] $G$ has a non-inner automorphism of order $p$;
			
			\item[\rm(2)] Exists a special subgroup $M$ of $G$ such that $|M|< |N|$;
			
			\item[\rm(3)] Exists a special subgroup $M$ of $G$ such that $|M|=|N|,|\mathrm{I}(C_{G}(N))|<|\mathrm{I}(C_{G}(M))|$. 
		\end{enumerate}
	\end{thm}
	\begin{pf}
		From Lemma \ref{kl}, there exists a normal subgroup $N_{1}$ of $G$ such that $\Omega_{1}(Z(N))Z(G)<N_{1}<N,|N/N_{1}|=p$ and $C_{G}(N)N/N_{1}$ non-cyclic.
		
		Since $C_{G}(N_{1})N_{1}/N_{1}$ is non-cyclic, $G$ has a non-inner automorphism of order $p$ from Lemma \ref{qk}.
			
		The proof is completed.
	\end{pf}
	
	\begin{lem}\label{xy}
		Let $N$ be a special subgroup of $G$ that $N/Z(G)\Omega_{1}(Z(N))$ is non-cyclic.  Suppose $\Omega_{1}(Z(N))$ is a exactly $n-G/N$-module where $n=d(Z(G))$. And all automorphisms induced by derivations in $Z^{1}(G/N,\Omega_{1}(Z(N)))$ are inner.  Assume $N_{1}$ is a normal subgroup of $G$ that $$ \Omega_{1}(Z(N))Z(G)\leq N_{1}< N=Z(N)N_{1}< C_{G}(N_{1})N.$$ If $\mathrm{H}^{1}(G/N_{1},\Omega_{1}(Z(N)))\geq \mathrm{H}^{1}(G/N,\Omega_{1}(Z(N)))\oplus\mathbb{F}_{p}^{2}$, $G$ has a non-inner automorphism of order $p$.
	\end{lem}
	\begin{pf}
		Denote $W=\Omega_{1}(Z(N))$.
		From Lemma \ref{ll}, suppose $C_{G}(N_{1})N/N$ is cyclic. Then $|\mathrm{I}(C_{G}(N_{1}))/\mathrm{I}(Z(N))|\leq p$. When $|\mathrm{I}(C_{G}(N_{1}))/\mathrm{I}(Z(N))|=1$, set $$\tau\in Z^{2}(G/N_{1},W)\backslash Z^{2}(G/N,W).$$ And $\phi(g)=g\tau(gN_{1})$. Then $\phi$ is an automorphism of order $p$ fixes $N_{1}$ not $N$. If $\phi$ is inner, there exists $h\in C_{G}(N_{1})$ such that $\phi(g)=g^{h}$ for all $g\in G$. Since $\phi$ fixes $N_{1}$ not $N$ and $ N=N_{1}Z(N)$, $ h\notin N$. Then $|\mathrm{I}(C_{G}(N_{1}))/\mathrm{I}(Z(N))|>1$. It's a contradiction.
		
		Suppose  $|\mathrm{I}(C_{G}(N_{1}))/\mathrm{I}(Z(N))|= p$. Set \begin{align*}
			\psi:\mathrm{I}(C_{G}(N_{1}))&\rightarrow Z^{1}(G/N_{1},W)\\
			x&\mapsto \d_{x} 
		\end{align*}
		where $\d_{x}(g)=g^{-1}g^{x}$. $\psi$ induce a monomorphism from $\mathrm{I}(C_{G}(N_{1}))/WZ(G)$ to $\mathrm{H}^{1}(G/N_{1},W)$. Since $\mathrm{H}^{1}(G/N_{1},\Omega_{1}(Z(N)))\geq \mathrm{H}^{1}(G/N,\Omega_{1}(Z(N)))\oplus\mathbb{F}_{p}^{2} $ and $$|\mathrm{I}(C_{G}(N_{1}))/\mathrm{I}(Z(N))|= p,$$ $ \psi$ is non-surjective. Set $\d\in Z^{1}(G/N_{1},W)\backslash \mathrm{Im}\psi$. And $\xi(g)=g\d(gN_{1})$.  $\xi$ is a non-inner automorphism of order $p$.
		
		The proof is completed.
	\end{pf}

	\begin{thm}\label{9.2} 
		Let $N$ be a special subgroup of $G$ that $N/Z(G)\Omega_{1}(Z(N))$ is non-cyclic.  Suppose $\Omega_{1}(Z(N))$ is a exactly $n-G/N$-module where $n=d(Z(G))$. And all automorphisms induced by derivations in $Z^{1}(G/N,\Omega_{1}(Z(N)))$ are inner. Suppose $ C_{G}(N)N/N$ is trivial. Set $N_{1}$ is a normal subgroup of $G$ that $\Omega_{1}(Z(N))Z(G)< N_{1}< N,N/N_{1}\cong\mathrm{C}_{p}$. If $C_{G}(N_{1})N_{1}/N_{1}$ is non-cyclic, then  
		
	\begin{enumerate}
		\item[\rm(1)] $G$ has a non-inner automorphism of order $p$;
		
		\item[\rm(2)] Exists a special subgroup $M$ of $G$ such that $|M|< |N|$;
		
		\item[\rm(3)] Exists a special subgroup $M$ of $G$ such that $|M|=|N|,|\mathrm{I}(C_{G}(N))|<|\mathrm{I}(C_{G}(M))|$. 
	\end{enumerate}
	\end{thm}
	\begin{pf}
		Denote $W=\Omega_{1}(Z(N))$.
		From Lemma \ref{ll}, suppose $C_{G}(N_{1})N/N$ is cyclic.   Since $C_{G}(N_{1})N_{1}/N_{1}$ is non-cyclic, $N< C_{G}(N_{1})N_{1}$.
		
		Since $N< C_{G}(N_{1})N_{1}$, there exists a normal subgroup $A$ of $G$ such that $N\leq A\leq C_{G}(N_{1})N_{1}$ and $A/N_{1}\cong\mathrm{C}_{p}^{2}$. Then $\mathrm{I}(C_{G}(N_{1}))\leq A$. 
		
		Claim $[x,g]=1 $ for all $x\in \Omega_{1}(Z(N)),g\in A $. Otherwise, there exists $x\in \Omega_{1}(Z(N))$ such that $[x,g]\neq 1$ for some $g\in A\backslash N$. Notice $ \Omega_{1}(Z(N))\leq N_{1}$. Thus $[y,g]=1$ for all $y\in \Omega_{1}(Z(N)),g\in C_{G}(N_{1})$. It's contradictory to $x\in \Omega_{1}(Z(N))$. Then $W=C_{W}(A)$.

		When $ n\geq 2$, claim $\mathrm{H}^{1}(G/N_{1},W)\geq \mathrm{H}^{1}(G/N,W)\oplus\mathbb{F}_{p}^{2}$. If $\mathrm{H}^{1}(G/N,W)\cong\mathrm{H}^{1}(G/A,W)$, from (2) in Theorem \ref{qq}, $\mathrm{H}^{1}(G/N_{1},W)\geq \mathrm{H}^{1}(G/N,W)\oplus\mathbb{F}_{p}^{2}$. If  $\mathrm{H}^{1}(G/A,W)\cong\mathbb{F}_{p}^{m} $ where $m<n$, from (1) in Theorem \ref{qq}, $\mathrm{H}^{1}(G/N_{1},W)\geq \mathrm{H}^{1}(G/N,W)\oplus\mathbb{F}_{p}^{2}$. The assertion is proved. Since $N=Z(N)N_{1}$, from Lemma \ref{xy}, $G$ has a non-inner automorphism of order $p$.
		
		When $n=1$, from Theorem \ref{dp} and Corllary \ref{yy} and $W=C_{W}(A)$, $$\mathrm{H}^{1}(G/N_{1},W)\cong\mathrm{H}^{1}(G/N,W)\oplus \mathrm{F}_{p}.$$ Set $\tau\in  Z^{1}(G/N_{1},W)\backslash Z^{1}(G/N,W)$. And $\phi(g)=g\tau(gN_{1})$ for all $g\in G$. Then $\phi$ is order $p$. If $\phi$ is inner ,there exists $ a\in A\backslash N$ such that $\phi(g)=g^{a}$ for all $g\in G$.  
Set $B=\langle a,N_{1}\rangle$. Then $B$ is a normal subgroup.

When $I(C_{G}(N))\leq N_{1}$, from Lemma \ref{qp}, either $B$ is a special subgroup of $G$ which has properties in (2) above or $G$ has a non-inner automorphism of order $p$.

When $\mathrm{I}(C_{G}(N))\nleqslant N_{1}$, since $\mathrm{I}(C_{G}(N))/\Omega_{1}(Z(N))Z(G)$ is cyclic, $N\neq \mathrm{I}(C_{G}(N))$. 
Then there exists a normal subgroup $N_{2}$ of $G$  
such that $\mathrm{I}(C_{G}(N))\leq N_{2}\leq N,|N/N_{2}|=p$. If $ C_{G}(N_{2})N_{2}/N_{2}$ is non-cyclic, we complete the proof from 
the above discussion. If $ C_{G}(N_{2})N_{2}/N_{2}$ is cyclic, we will discuss in Theorem \ref{tt}.
		
		The proof is completed.
	\end{pf}
	
	\begin{thm}\label{tt}Let $N$ be a special subgroup of $G$ that $N/Z(G)\Omega_{1}(Z(N))$ is non-cyclic.  Suppose $\Omega_{1}(Z(N))$ is a exactly $n-G/N$-module where $n=d(Z(G))$. And automorphisms induced by derivations in $Z^{1}(G/N,\Omega_{1}(Z(N)))$ are inner. Suppose $ C_{G}(N)N/N$ is trivial.
		Set $N_{1}$ is a normal subgroup of $G$ that $ \Omega_{1}(Z(N))Z(G)\leq N_{1}\leq N,|N/N_{1}|=p$. If $C_{G}(N_{1})N_{1}/N_{1}$ is non-trivial cyclic $p$-group,  then one of following holds:
		\begin{enumerate}
			\item[\rm(1)] $G$ has a non-inner automorphism of order $p$;
			
			\item[\rm(2)] Exists a special subgroup $M$ of $G$ such that $|M|< |N|$;
			
			\item[\rm(3)] Exists a special subgroup $M$ of $G$ such that $|M|=|N|,|\mathrm{I}(C_{G}(N))|<|\mathrm{I}(C_{G}(M))|$. 
		\end{enumerate}
	\end{thm}
	\begin{pf}
		Since $C_{G}(N_{1})N_{1}/N_{1}$ is cyclic, there are two cases:
		\begin{align*}
			&(1)\;N\leq C_{G}(N_{1})N_{1},\\
			&(2)\;C_{G}(N_{1})N_{1}\cap N=N_{1}.\\
		\end{align*}
		Denote $W=\Omega_{1}(Z(N))$.
		
		Firstly, discuss $C_{G}(N_{1})N_{1}\cap N=N_{1}$. When $\mathrm{I}(C_{G}(N_{1}))\leq N_{1}$ $N_{1}$ has proerties in {\rm(2)} in the above. Next, we discuss the case $ \mathrm{I}(C_{G}(N_{1}))\nleqslant N_{1}$. Set $ b\in \mathrm{I}(C_{G}(N_{1}))\backslash N$ that $b^{p}\in N_{1}$. And $B=\langle b,N_{1}\rangle$.
		 Then $\mathrm{I}(C_{G}(B))\leq B$ and $\mathrm{I}(C_{G}(N))$ is a proper subgroup of  $\mathrm{I}(C_{G}(B))$. From Lemma \ref{qp}, either $G$ has a non-inner automorphism of order $p$ or $B$ has properties in {\rm(3)} above.

		Secondly, discuss $C_{G}(N_{1})N_{1}\cap N=N$. Claim $C_{G}(N_{1})N_{1}=N$. Otherwise, set $a\in C_{G}(N_{1})\backslash N$ that $a^{p}\in N$. Since $ C_{G}(N_{1})N_{1}/N_{1}$ is cyclic, $N=\langle a^{p},N_{1}\rangle$. Thus $a\in C_{G}(N)$. It's contradictory to $C_{G}(N)\leq N$. When $\mathrm{H}^{1}(G/N_{1},W)\geq \mathrm{H}^{1}(G/N,W)\oplus\mathbb{F}_{p}$, there exists $\tau\in Z^{1}(G/N_{1},W)\geq Z^{1}(G/N,W)$. Set $\psi$ is automorphism induced by $\tau$. Since $C_{G}(N_{1})$ is abelian, $\psi$ is non-inner. 
		
		Suppose $\mathrm{H}^{1}(G/N_{1},W)\cong \mathrm{H}^{1}(G/N,W)$. 
		When $ N$ is abelian, if $n\geq2$, there exists a normal subgroup $N_{2}$ of $G$ such that $WZ(G)\leq N_{2}\leq N_{1}$ and $N/N_{2}\cong\mathrm{C}_{p}^{2}$. From (2) in Theorem \ref{qq}, $\mathrm{H}^{1}(G/N_{2},W)\geq \mathrm{H}^{1}(G/N,W)\oplus \mathbb{F}_{p}^{2}$. By Lemma \ref{xy}, $G$ has a non-inner automorphism of order $p$. 
		
		If $n=1$, since $N/WZ(G)$ is non-cyclic, there exists a normal subgroup $T$ of $G$ such that $\mathrm{I}(C_{G}(N))\leq T\leq N,|N/T|=p$. 
		When $C_{G}(T)T/T$ is non-cyclic, from the proof of Theorem \ref{9.2}, we got the conclusion.
		When $C_{G}(T)T/T$ is cyclic, from the above discussion, we assume $C_{G}(T)=N$ and $\mathrm{H}^{1}(G/T,W)\cong \mathrm{H}^{1}(G/N,W)$. Then $I(C_{G}(T))=I(C_{G}(N))$. 
		$T$ is a special subgroup that satisfying (2). 
		
		When $N$ is non-abelian, there exists a normal subgroup $N_{2}$ of $G$ such that $Z(N)\leq N_{2}\leq N,|N/N_{2}|=p$. If $C_{G}(N_{2})N_{2}/N_{2}$ is non-cyclic, from Theorem \ref{9.2}, the result can be got. 
		Suppose $C_{G}(N_{2})N_{2}/N_{2}$ is non-trivial cyclic. From first case, the result can be got.
		
		The proof is completed.
	\end{pf}
	
	Before discussing the last case, we must provide some lemmas.
	
	\begin{lem}\label{cct}
		$N$ is a normal subgroup of $G$ that satisfying $ N\leq Z(G),N\cong C_{p}^{2}$ and $G/N$ non-cyclic. Then $$\tau(g)\notin  J(\mathbb{F}_{p}(G/N))$$ for any $\tau\in Z^{1}(G,\mathbb{F}_{p}(G/N))\backslash Z^{1}(G/N,\mathbb{F}_{p}(G/N))$ and $g\in G$.
		\end{lem}
	\begin{pf}
	We prove this by contradiction. Suppose there exists  $$\tau\in Z^{1}(G,\mathbb{F}_{p}(G/N))\backslash Z^{1}(G/N,\mathbb{F}_{p}(G/N))$$ and $g\in G\in N$ such that $ \tau(g)\in \mathbb{F}_{p}(G/N)\backslash J(\mathbb{F}_{p}(G/N))$. Then $\tau(g^{p^{s}})=\tau(g)(\overline{g}-1)^{p^{s}-1}$ where $\overline{g}=gN$ and $o(\overline{g})=p^{s}$. Since $G/N$ is non-cyclic, $(\overline{g}-1)^{p^{s}-1}\notin \mathbb{F}_{p}(G/N)^{G}$. Then $ \tau(g^{p^{s}})\notin \mathbb{F}_{p}(G/N)^{G}$. It's a contradiction.
	\end{pf}
	
	\begin{lem}\label{ook}
		Let $G=\langle a,b|a^{p^{n}}=b^{p}=c^{p}=1,c=[a,b],[c,a]=[c,b]=1\rangle$ where $n\geq 2$. And $N=\langle a^{p^{s}},b,c\rangle$ where
		$s\leq n-1$ and $p^{s}\geq 3$. Then $\mathrm{H}^{1}(G,J(\mathbb{F}_{p}(G/N)))\cong \mathrm{H}^{1}(G/\langle c\rangle,J(\mathbb{F}_{p}(G/N)))\oplus\mathbb{F}_{p}$.
	\end{lem}
	\begin{pf}
		Set $x\notin \mathbb{F}_{p}(G/N)^{G/N}$ that satisfying $ x^{g-1}\in \mathbb{F}_{p}(G/N)^{G/N}$ for all $g\in G$.
		Since $|G/N|= p^{s}\geq 3$, $x\in J(\mathbb{F}_{p}(G/N))$.
		Set $ \tau(1),\tau(b)=x$. It's easy to verify  that $\tau$ can be lifted to a dervation  from $G$ to $J(\mathbb{F}_{p}(G/N))$. And $\tau(c)=x^{1-a}\neq 0$.
		
		The proof is completed.
	\end{pf}

\begin{lem}\label{klp}
		Let $G$ be one of following groups
		$$D_{2^{n}},Q_{2^{n}},SD_{2^{n}}$$ where $n\geq 4$. 
		$M$ is a $1-G/T$-module that $\mathrm{H}^{1}(G,M)\cong \mathrm{H}^{1}(G/T,M)$ where $T$ is a minimal normal subgroup of $G$. And $G'$ acts on $M$ non-trivially.  Then 
		
		\begin{enumerate}
			\item[\rm(1)]  When $G$ is isomorphic to one of following groups :$ D_{2^{n}},SD_{2^{n}}$, $d(_{\tau}G.M)=d(G)+1$ for any $\tau\in Z^{2}(G,M)$
			\item[\rm(2)] 
			when $G\cong Q_{2^{n}}$ and $dim_{\mathbb{F}_{2}} M^{\Phi(G)}=2$, $M$ is $L$-module isomorphic to a quotietnt module of $\mathbb{F}_{2}(L)$ that $dim_{\mathbb{F}_{2}}M\geq 3$ where $L$ is a cyclic maximal subgroup of $G$.
		\end{enumerate}
\end{lem}
\begin{pf}
	From Lemma \ref{kj}, $\mathrm{H}^{1}(G,M)\cong \mathbb{F}_{2}$. 
	When $\tau\in B^{1}(G,M)$, $_{\tau}G.M\cong M\rtimes G$. Then $ d(_{\tau}G.M)=d(G)+1$.
	
	Next, we only discuss the case $ \tau\notin  B^{1}(G,M)$.
	Since $\mathrm{H}^{1}(G,M)\cong \mathbb{F}_{2}$, we only prove $d(_{\tau}G.M)=d(G)+1$ for any $ \tau \in Z^{2}(G,M)\backslash B^{2}(G,M)$.
	
	When $G$ is isomorphic to $D_{2^{n}}$, notice $D_{2^{n}}$ is a quotient group of $D_{2^{n}}$. Set $ K=(M\rtimes D_{2^{n+1}})/\langle ab \rangle $
 where $a\in Z(D_{2^{n+}})\backslash \{1\},b\in M^{G}\backslash \{1\}$. Notice $K$ is a group that $M$ extension by $G$. Since $ T$ acts on $M$ trivially, the extension is non-trivial.

When $G$ is isomorphic to $SD_{2^{n}}$, we need prove $ M^{G'}\leq J_{G}(M)$.  Set $y\in G$ that $\langle y\rangle$ is a cyclic maximal subgroup. Since $G'$ acts on $M$ non-trivially, there exist $x\in M $ such that $ x^{y^{2}-1}\in M^{G}\backslash \{1\}$ and $x^{g-1}\in M^{\Phi(G)}$ for all $g\in G$.  Then $dim_{\mathbb{F}_{2}}M\geq 3$ and $(x^{y-1})^{y-1}\neq 0$. Set $M_{1}=\langle x\rangle\oplus M^{G}$. Then $M_{1}$ is a $G$-submodule that $dim_{\mathbb{F}_{2}}M_{1}=2$.  From Theorem 2 in\cite{Nin}, set $$L\cong (M\rtimes G_{20})/\langle b^{2}x, a^{2^{m-3}}u\rangle$$ where $u\in M^{G}\backslash \{1\}$. Then $L$ is a group extension that $M$ by $G$. Notice $ d(L)=3$. Since $ T$ acts on $M$ trivially, the extension is non-trivial.

Next, we suppose that  $G$ is isomorphic to $Q_{2^{n}}$. It's similar to prove $ dim_{\mathbb{F}_{2}}M\geq 3$. Notice  $x^{a^{2}-1}\neq 0$  for any $x\in M\backslash M^{\Phi(G)}$ where $\langle a\rangle$ is a cyclic maximal subgroup of $G$. Then $x^{a-1}\neq 0$. Since $dim_{\mathbb{F}_{2}} M^{\Phi(G)}=2$, $M$ is a $\langle a\rangle$-module isomorphic to a quotietnt module of $\mathbb{F}_{2}(L)$.

The proof is completed.
\end{pf}

\begin{lem}\label{iio}
	Let $G\cong G_{6}$ where $G_{6} $ defined in Theorem 2 in \cite{Nin}. $H$ is finite $p$-group that has a minimal normal subgroup $N$ satisfying $H/N\cong G$. If the minimal index of  cyclic subgroups of $H$ is $8$, then $H$ is isomorphic to one of following groups:
	\begin{enumerate}
		\item[\rm(1)] $C_{2}\times G_{6}$;
		\item[\rm(2)] $\langle a,b|a^{2^{n-3}}=b^{8}=1, a^{b} =a^{-1}b^{4q}\rangle$;
		\item[\rm(3)] $ \langle a,b|a^{2^{n-3}}=b^{4}=c^{2}=1, a^{b}=a^{-1}c,[c,a]=[c,b]=1\rangle$
	\end{enumerate}
where $ n\geq 6,q=0,1$.
\end{lem}
\begin{pf}
Notice the groups in the above have quotient group which is isomorphic to $G_{6}$.	We only need to prove that $H$ is isomorphic to one of the above groups. When $d(H)=3$, $H\cong C_{2}\times G_{6}$. Next, we suppose $d(H)=2$. Set $$G=\langle x,y|x^{2^{n-2}}=y^{4}=1,x^{y}=x^{-1}.\rangle$$   
Set $\pi$ is surjective group homomorphism from $H$ to $G$. And $a,b\in H$ that $\pi(a)=x,\pi(b)=y$. Since the minimal index of  cyclic subgroups of $H$ is $8$, $a^{2^{n-2}}=1$. When $b^{4}\neq 1$, $a^{b}=a^{-1}b^{4q}$ where $q=0,1$. Their groups are described in  the above (2). When $b^{4}=1$, there is $z\in ker\pi$. Then $a^{b}=a^{-1}z^{s}$ where $s=0,1$. If $s=0$, $z\notin \langle a,b\rangle$. It's a contradiction. Then $$H=\langle a,b|a^{2^{n-2}}=b^{4}=z^{2}=1,a^{b}=a^{-1}z,[z,a]=[z,b]=1\rangle.$$	
	The proof is completed.
\end{pf}
	
	\begin{thm}\label{uui}
		Let $G$ be a finite $p$-group with unique elementary abelian normal $p$-subgroup. $N$ is the minimal normal subgroup of $G$. $M$ is a $1-G/N$-module that satisfying $$\mathrm{H}^{1}(G,M)\cong \mathrm{H}^{1}(G/N,M).$$
		Then $_{\tau}G.M$ has a non-inner automorphism of order $p$ for any $\tau\in Z^{2}(G,M)$ that satisfying $d(Z(_{\tau}G.M))=1$.
	\end{thm}
\begin{pf}
Since $G$ is a finite $p$-group with unique elementary abelian normal $p$-subgroup, from chapter 5 Theorem 4.10 in \cite{gor},  $G$ is isomorphic to one of the following groups:
$$ D_{2^{m}},Q_{2^{n}},SD_{2^{m}}$$ where $n\geq 3,m\geq 4$. 

When $G'$ acts on $M$ trivially, $\Phi(_{\tau}G.M)$  is abelian. If $d(\Phi(_{\tau}G.M))=3$, it's easy to prove $G$ has non-inner auomorphism of order $p$. If $d(\Phi(_{\tau}G.M))=2$, from \cite{mohsen2}, $G$ has non-inner automorphism of order $2$.
Then we suppose $G'$ acts on $M$ non-trivially in the following discussion. Since $N$ acts on $M$ trivially, $G$ is isomorphic to one of the following groups:
$$ D_{2^{n}},Q_{2^{n}},SD_{2^{n}}$$ where $n\geq 4$.

From Lemma \ref{klp}, $d(_{\tau}G.M)=d(G)+1$ when $G$  is isomorphic to one of the following groups: $D_{2^{n}},SD_{2^{n}}$. Since $\mathrm{H}^{1}(_{\tau}G.M/J_{G}(M),M^{G})\cong \mathbb{F}_{p}^{3}$, there exists $\tau\in Z^{1}(G)(_{\tau}G.M/J_{G}(M),J_{G}(M))\backslash Z^{1}(G)(G,J_{G}(M))$ such that $\tau(g)\in M^{G}$ for all $g\in G$.
Set $\psi(g)=g\tau(gJ_{G}(M))$. Then $\psi$ is order $p$. It's easy to prove $\psi$ is non-inner.

Next, we discuss that $G$ is isomorphicto $Q_{2^{n}}$. Since $d(_{\tau}G.M)=d(G)+1$, it's similar to prove $_{\tau}G.M$ has a non-inner automorphism of order $2$. Then we suppose $d(_{\tau}G.M)=d(G)$.  
 Since there isn't finite $2$-group with cyclic maximal subgroup such that has proper quotient group which is isomorphic to $Q_{2^{n}}$, the extension of $C_{2}$ by $Q_{2^{n}}$ has cyclic subgroup which has index $8$. From Theorem 2 in \cite{Nin}, $_{\tau}G.M/J_{G}(M)$ is isomorphic to one of following groups: $G_{2},G_{6},G_{14}$ where $G_{2},G_{6},G_{14}$ defined in \cite{Nin}. When $_{\tau}G.M/J_{G}(M)\cong G_{2}$, $d(_{\tau}G.M/J_{G}(M))=d(G)+1$. It's easy to prove $_{\tau}G.M$ has non-inner automorphism of order $p$ similar to above.
 
  we assert $_{\tau}G.M/J_{G}(M)\cong G_{14}$. Otherwise $_{\tau}G.M/J_{G}(M)\cong G_{6}$ Since $dim_{\mathbb{F}_{2}}M\geq 3$, there exists $G$- submodule $M_{1}$ of $M$ such that $ dim_{\mathbb{F}_{2}}M/M_{1}=2$ and $M^{G}\leq M_{1}$. Then the cyclic subgroups of $_{\tau}G.M/M_{1}$ with maximal order have index $p^{3}$. From Lemma \ref{iio}, $_{\tau}G.M/M_{1}$ is isomorphic to $$\langle a,b|a^{2^{n-3}}=b^{4}=c^{2}=1, a^{b}=a^{-1}c,[c,a]=[c,b]=1\rangle.  $$
 Since $ M$ is cyclic module, $G$ acts on $M/M_{1}$ non-trivially. Notice $ d(Z(G))=3$. It's contradictory to $G$ acts on $M/M_{1}$ non-trivially.

  Now, we begin to discuss the case $_{\tau}G.M/J_{G}(M)\cong G_{14}$, since $G'$ acts $M$ non-trivially, there exists $x\in M^{\Phi(G)}$ such that $x^{a-1}\neq 0$ where $ \langle a\rangle$ is a cyclic maximal subgroup of $G$. Notice $$G_{14}/\langle a^{4}\rangle\cong \langle u,v|u^{4}=v^{2}=t^{2}=1,t=[u,v],[t,u]=[t,v]=1\rangle.$$
Set $\d(v)=x,\d(u)=0$. It's easy to verfy $\d$ can be extent a dervation from $G_{14}/\langle a^{4} $ to $ M^{\Phi(G)}$. Then there exists $\kappa\in Z^{1}(_{\tau}G.M/J_{G}(M),M^{\Phi(G)})\backslash Z^{1}(G,M^{\Phi(G)})$. Notice $M^{Phi(G)}\leq J_{G}(M)$. Then $\kappa\in Z^{1}(_{\tau}G.M/J_{G}(M),J_{G}(M))\backslash Z^{1}(G,J_{G}(M))$. Set $\kappa_{1}(g)=g\kappa(gJ_{G}(M))$. Then $\kappa_{1}$ is automorphism of order $2$. Since $ G$ has unique elementary abelian group, $\kappa_{1}$ is non-inner.

The proof is completed.
\end{pf}

\begin{lem}\label{lop}
	Let $G$ be a finite $p$-group with non-cyclic elementary abelian normal $p$-group. $N$ is a cyclic normal subgroup of $G$. $Q$ is a $1-G$-module that satisfying $\mathrm{H}^{1}(G,Q)\cong \mathrm{H}^{1}(G/N,Q),|Q|=p^{2}$.Then $_{\tau}G.Q$ has non-inner automorphism of order $p$ for any $\tau\in Z^{1}(G,Q)$.
\end{lem}
	\begin{pf}
		From Lemma \ref{kko},  there exists maximal subgroup  $T$ of $G$ such that $T$ acts on $M$ trivially. Then $C_{G}(Z(\Phi(G)))\neq \Phi(G)$. From main theorem in \cite{marian}, $G$ has a non-inner automorphism of order $p$.
	\end{pf}

	\begin{thm}\label{xpl}
		Let $N$ be a special subgroup of $G$ that $N/Z(G)\Omega_{1}(Z(N))$ is cyclic. Suppose $\Omega_{1}(Z(N))$ is a exactly $n-G/N$-module that $d(Z(G))=n$. And all automorphisms induced by derivations in $Z^{1}(G/N,\Omega_{1}(Z(N)))$ are inner. Then one of following holds:
		\begin{enumerate}
		\item[\rm(1)] $G$ has a non-inner automorphism of order $p$;
		
		\item[\rm(2)] Exists a special subgroup $M$ of $G$ such that $|M|< |N|$;
		
		\item[\rm(3)] Exists a special subgroup $M$ of $G$ such that $|M|=|N|,|\mathrm{I}(C_{G}(N))|<|\mathrm{I}(C_{G}(M))|$. 
		\end{enumerate}
	\end{thm}
	\begin{pf}
		Notice $\Omega_{1}(Z(N))$ is exactly $n-G/N$-module. From $N/Z(G)\Omega_{1}(Z(N))$ is cyclic, $N$  is abelian. And $n=1$. Denote $W=\Omega_{1}(Z(N))$.
		
		When $C_{G}(N)N/N$ is non-trivially cyclic, there exists a normal subgroup $M$ of $G$ such that $\Omega_{1}(Z(N))Z(G)\leq M<N$ such that $|N/M|=p$. If $C_{G}(N)N/M$ is non-cyclic, from Lemma \ref{qk}, $G$ has a non-inner automorphism of order $p$.
		
		If $C_{G}(N)N/M$ is cyclic,  then $C_{G}(N)N/W$ is cyclic.
		Thus $$\mathrm{H}^{1}(G/W,W)\cong\mathrm{H}^{1}(G/C_{G}(N)N,W)\cong\mathbb{F}_{p}.$$
		From Lemma \ref{jx}, $G/W$ is non-abelian.  And $W$ is cyclic $G$ module by Lemma \ref{kj}.
	
	When $G/W$ has unique normal elementary abelian subgroup, from Theorem \ref{uui}, $G$  has a non-inner automorphism of order $p$. 
	
	When $G/W$ has normal elementary abelian subgroup with order $p^{2}$, If $|W|=p^{2}$, from Lemma \ref{lop}, $G$ has a non-inner automorphism of order $p$.
	
	Next, we discuss the case $|M|\geq p^{3}$.
		From Theorem \ref{du}, $\mathrm{H}^{1}(G/\mathrm{J}_{G}(W),\mathrm{J}_{G}(W))\cong\mathrm{H}^{1}(G/W,\mathrm{J}_{G}(W))\oplus\mathbb{F}_{p}$.   Set $\d\in Z^{1}(G/\mathrm{J}_{G}(W),\mathrm{J}_{G}(W))\backslash Z^{1}(G/W,\mathrm{J}_{G}(W))$ and $\phi(g)=g\d(g\mathrm{J}_{G}(W))$. $\phi$ is order $p$. If $\phi$ is inner, there exists $x\in C_{G}(\mathrm{J}_{G}(W))\backslash C_{G}(N)N$ such that $\phi(g)=g^{x}$ for all $g\in G$. Set $B=\langle x,N \rangle$. $B$ is normal subgroup of $G$ that $B/M\cong\mathrm{C}_{p}\times C_{p^{s}}$ where $s\geq 1$. Since $[x,g]\in \mathrm{J}_{G}(W)$ for all $g\in G$, $xW\leq \Omega_{1}(Z(G/W))$. From $G/W$ is non-abelian and $C_{G}(N)N\leq\Phi(G)$, $ G/B$ is non-cyclic. From Theorem \ref{dp}, $ |C_{W}(B)|\geq p^{p}$. From Corllary \ref{rty} and $\mathrm{H}^{1}(G/M,W)\cong\mathrm{H}^{1}(G/N,W)$, $ |W|=|C_{W}(B)|^{p}\geq p^{p^{2}}$. Then $$C_{W}(B)< \mathrm{J}_{G}(W).$$ 
		Thus there exists $y\in \mathrm{J}_{G}(W)\backslash C_{W}(B)$ such that $[x,y]\neq 1$. It's contradictory to $\d\in Z^{1}(G/\mathrm{J}_{G}(W),\mathrm{J}_{G}(W))\backslash Z^{1}(G/W,\mathrm{J}_{G}(W))$. Then $\phi$ is non-inner.
		
		When $ C_{G}(N)N/N$ is trivial,  there exists a normal subgroup $T$ of $G$ such that $\Omega_{1}(Z(N))Z(G)\leq T<N$ such that $|N/T|=p$. if $\mathrm{I}(C_{G}(N))\leq T$, then $M$ has properties in the above (2). Next, suppose $ N=\mathrm{I}(C_{G}(N))$. Then $T=WZ(G)$. 
		
		When $\mathrm{H}^{1}(G/T,W)\cong \mathrm{H}^{1}(G/N,W)\oplus \mathbb{F}_{p}$, thre exists $\kappa\in Z^{1}(G/T,W)\backslash Z^{1}(G/N,W)$. Set $\kappa_{1}(g)=g\kappa(gT)$. Then $\kappa_{1}$ is order $p$. When $\kappa_{1}$ is inner, there exists $f\in G\backslash N$ such that $\kappa_{1}(g)=g^{f}$. Set $L=\langle N,f\rangle$. We assert $C_{W}(L)=W$. Otherwise, it's contradictiory to $W\leq T$ and $f\in C_{G}(T)$.
		Notice $N=I(C_{G}(N))$ and $C_{G}(N)=N$. Then $\mathrm{H}^{1}(G/N,W)\oplus \mathrm{H}^{1}(G/L,W)\oplus \mathbb{F}_{p}$. Thus $ \mathrm{H}^{1}(G/L,W)=0$. This imply $W\cong\mathbb{F}_{p}(G/L)$. When $G/L$ is non-cyclic, from Lemma \ref{cct}, $L/J_{G}(W)\leq Z(G/J_{G}(W))$. From Lemma \ref{cct}, there exists $\a\in Z^{1}(G/J_{G}(W),J_{G}(W))\backslash Z^{1}(G/W,J_{G}(W))$.  Set $\a_{1}(g)=g\a(gJ_{G}(W))$. $\a_{1}$ is order $p$. We assert $\a_{1}$ is non-inner. Otherwise there exists $q\in G\L$ such that $g^{q}=\a_{1}(g)$. Set $L_{1}=\langle L,q\rangle$. Since $W\cong \mathbb{F}_{p}(G/L)$, $|W|=|C_{W}(L_{1})|^{p}$. From $G/L$ is non-cyclic, $C_{W}(L_{1})<J_{G}(W)$. It's contradictory to $q\in C_{G}(W)$.
		
		When $G/L$ is cyclic, $f\notin \Phi(G)$.  Then $G/L$ is abelian.   Notice $G/L$ acts on $L/T$ trivially. Then ether $ G/T\cong C_{p^{s}}\times C_{p}\times C_{p}$ or $G/T\cong C_{p^{t}}\times C_{p}$ where $s\geq 1,t\geq 2$. Since $N\leq \Phi(G)$, $G/T\cong C_{p^{t}}\times C_{p}$. Since $d(N/W)=1$ and $[f,g]\in W$ for all $g\in G$, $G/W$ is abelian. From Lemma \ref{jx}, $d(G)=d(G/N)+1$. It's contradictory to $N\leq \Phi(G)$. The case that $d(G/L)=1$ isn't valid.		
		
		From discussion in the above  we suppose that $\mathrm{H}^{1}(G/N,W)\cong \mathrm{H}^{1}(G/W,W),|W|\geq p^{3}$ and $G/W$ has  normal elementary abelian subgroup with order $p^{2}$. Similar to discussion in the above, there exists a automorphism of  order $p$ $\psi$ such that $\psi$ fix $J_{G}(W)$ not $W$.
		If $\psi$ is inner, there exists $y\in G\backslash N$ such that $\psi(g)=g^{y},y^{p}\in Z(G)$ for all $g\in G$. Set $ D=\langle y,N\rangle$. Since $G/W$ is non-abelian, $G/D$ is non-cyclic. Since $\mathrm{H}^{1}(G/T,W)\cong \mathrm{H}^{1}(G/N,W)$, $|W|=|C_{W}(D)|^{p}$ from Corllary \ref{rty}. Since $G/D$ is non-cyclic, $C_{W}(D)\geq p^{2}$. Then $C_{W}(D)< J_{G}(W)$. There exists $h\in J_{G}(W)\backslash C_{W}(D)$. Then $[h,y]\neq 1$. It's a contradiction. $\psi$ is non-inner.
		
		The proof is completed.
	\end{pf}
	
	Summarize the above conclusions, the following Theorem can be got.
	
	\begin{thm}\label{ui}
		Suppose $N$ is a special subgroup of $G$. Then one of following holds:
		
		\begin{enumerate}
			\item[\rm(1)] $G$ has a non-inner automorphism of order $p$;
			
			\item[\rm(2)] Exists a special subgroup $M$ of $G$ such that $|M|< |N|$;
			
			\item[\rm(3)] Exists a special subgroup $M$ of $G$ such that $|M|=|N|,|\mathrm{I}(C_{G}(N))|<|\mathrm{I}(C_{G}(M))|$. 
		\end{enumerate}
	\end{thm}
	\begin{pf}
		Assume $n=d(Z(G))$. Denote $W=\Omega_{1}(Z(N))$.	From Theorem \ref{5.5}, suppose $\Omega_{1}(Z(N))$ is $ n-G/N$-module. If $ \mathrm{H}^{1}(G/N,W)\cong\mathbb{F}_{p}^{m}$ where $m<n$, by Theorem \ref{hh}, either $G$ has a non-inner automorphism of order $p$ or exists a special subgroup $N_{1} $ of $G$ such that $N_{1}<N$ when $C_{G}(N)\leq N$. When $C_{G}(N)\nleqslant N$, from Theorem \ref{gg} and Lemma \ref{ty} , $G$ has a non-inner automorphism of order $p$.
		
		If $\mathrm{H}^{1}(G/N,W)\cong\mathbb{F}_{p}^{n}$, from Theorem \ref{xx}, Theorem \ref{9.2} ,Theorem \ref{tt} and Theorem \ref{xpl}, one of following holds:
		
		{\rm(1)}\; $G$ has a non-inner automorphism of order $p$;
		
		{\rm(2)}\; Exists a special subgroup $M$ of $G$ such that $|M|< |N|$;
		
		{\rm(3)}\; Exists a special subgroup $M$ of $G$ such that $|M|=|N|,|\mathrm{I}(C_{G}(N))|<|\mathrm{I}(C_{G}(M))|$.
		
		The proof is completed.
	\end{pf}
	
	\begin{cor}\label{18}
		Any finite $p$-group of order at least $p^{2}$ has a non-inner automorphism of order $p$.
	\end{cor}
	\begin{pf} When $G$ is abelian, from Gasch{\"u}tz Wolfgang's classic proposition in \cite{gas}, $G$ has a outer automorphism of order $p$.  Since $Aut(G)\cong Aut(G)/Inn(G)$, $G$ has a non-inner automorphism of order $p$.
		
		Next, suppose $G$ is non-abelian. From main result in \cite{marian},
		suppose $G$ is a non-abelian finite $p$-group that $C_{G}(\Phi(G))\leq\Phi(G)$.
		Denote $S_{G}$ as the set consists of special subgroup of $G$ which has minimal order. Set $A\in S_{G}$ that $|\mathrm{I}(C_{G}(A))|\geq |\mathrm{I}(C_{G}(N))|$ for all $N\in S_{G}.$ From Theorem \ref{ui} and the selection of $A$, $G$ has a non-inner automorphism of order $p$. 
		
		The proof is completed.
	\end{pf}

\end{document}